\documentclass[twocolumn]{autart}    

\usepackage{graphicx}          
 \usepackage{mathrsfs}
\usepackage{amsfonts}
\usepackage{amssymb}
\usepackage{amsmath}
\usepackage{subfigure}
\usepackage{color}
\usepackage{float}
\usepackage{wrapfig}
\usepackage{soul}
\def\disp{\displaystyle}

\def\dref#1{(\ref{#1})}
\def\crr{\cr\noalign{\vskip0mm}}

\usepackage[figuresright]{rotating}

\newtheorem{theorem}{Theorem}[section]                   

\newtheorem{lemma}{Lemma}[section]
\newtheorem{remark}{Remark}[section]

\usepackage{booktabs}

\begin{document}

\begin{frontmatter}
\title{Event-Triggered Active Disturbance Rejection Control for Uncertain Random Nonlinear Systems 
} 

%

\author[foshan]{Ze-Hao Wu}\ead{zehaowu@amss.ac.cn},
\author[Guangzhou]{Feiqi Deng}\ead{aufqdeng@scut.edu.cn},    
\author[foshan]{Pengyu Zeng}\ead{pyzeng@outlook.com},             
\author[Changsha]{Hua-Cheng Zhou}\ead{hczhou@amss.ac.cn},  
\author[Chongqing]{Hongyi Li}\ead{lihongyi2009@swu.edu.cn}

\address[foshan]{School of Mathematics and Big Data, Foshan University, Foshan 528000, China}  
\address[Guangzhou]{Systems Engineering Institute, South China University of Technology, Guangzhou 510640, China}
\address[Changsha]{ School of Mathematics and Statistics, HNP-LAMA, Central South University, Changsha 410075, China}
\address[Chongqing]{School of Electronic and Information Engineering and Chongqing Key Laboratory of Generic Technology and System of Service Robots, Southwest University, Chongqing 400715, China}



\begin{abstract}                          
In this paper, event-triggered active disturbance rejection control (ADRC) is first addressed
 for a class of uncertain random nonlinear systems driven by bounded noise and colored noise.
The event-triggered extended state observer (ESO) and ADRC controller are designed, where two respective event-triggering mechanisms
 with a fixed positive lower bound for the inter-execution times are proposed. The random total disturbance representing the coupling of
nonlinear unmodeled dynamics, external deterministic disturbance, bounded noise, and colored noise is estimated in real
time by the  event-triggered ESO and compensated in the event-triggered feedback loop. Both the mean square and almost
surely practical convergence of the closed-loop systems is shown with rigorous theoretical analysis. Finally, some
numerical simulations are implemented to validate the proposed control scheme and theoretical results.
\end{abstract}

\end{frontmatter}

\section{Introduction}\label{section1}
\vskip-3mm
Disturbance rejection for uncertain systems has been attracting increasing attention in the past decades since disturbances
and uncertainties are widespread in engineering applications. Active disturbance rejection control (ADRC)
is a novel control technology based on estimation/compensation strategy, proposed by Han in the late 1980s \cite{Han2009}.
 The key component of ADRC scheme is extended state observer (ESO),
which aims at estimating in real time not only unmeasurable states but also total disturbance representing total effects of internal uncertainties and
external disturbances. Based on the estimates obtained via ESO, the ADRC controller composed of
a feedback controller and a compensator is devised to achieve disturbance rejection  and desired control objective for the plant,
where the total disturbance is promptly estimated and compensated in the closed-loop.
Since ADRC is proposed, it has been applied in extensive engineering technology applications like
 general purpose control chips manufactured by Texas Instruments \cite{Ti},
 DC-DC power converter \cite{DC-DC}, permanent magnet synchronous motor \cite{permanentmagnetsynchronousmotor}, Delta robot control \cite{robotcontorl},
 and power plant \cite{powerplant}, etc.

On the one hand, over the past twenty years much progress  on the theoretical foundation of ADRC for uncertain systems has been made including the convergence analysis of ESO and ADRC's closed-loop,
which can be found in the convergence analysis of ESO for a class of uncertain nonlinear systems in \cite{zhaoESO2011}, the convergence analysis of ADRC's closed-loop
of uncertain random systems in \cite{wu1}, and more all-round designs and theoretical analysis of ADRC for uncertain nonlinear systems in the monographs \cite{zhaoADRC2016,Sira2017} and the references therein.
On the other hand, networked control systems (NCSs), where sensors, actuators, and controllers are spatially distributed,
 transmit data over a communication network with superiority in reducing installation costs, obtaining higher dependability, and improving system flexibility, and then have been applied extensively.
It is important to note that  the bandwidth of the wireless sensor-actuator networks and computation resources are limited in NCSs.
Event-triggering mechanisms (ETMs) have been widely applied in improving communication efficiency for NCSs
 with ensuring control quality, see, e.g., \cite{ETM1,ETM2,ETM3,ETM4,ETM5} and the references therein.
The ESO designs via the continuously transmitted output measurement and the ADRC controller designs based on continuous-time estimates obtained by ESO
like those in aforementioned literature may become impracticable for uncertain systems in networked environment. This further promotes the
recent development of event-triggered ADRC for uncertain systems, see, e.g.,
\cite{eventADRC2017,eventADRC2019,eventADRC2021} and the references therein. In the framework of event-triggered control,
information transmission and control signal renovating occur only when essential for the system, leading to momentous
efficiency in saving communication/computation resources.

It is widely acknowledged that disturbances and uncertainties are more likely to be of random characteristic in engineering applications.
Many important developments have been made in the feedback control of stochastic systems driven by white noise described by stochastic differential equations (SDEs)
in the past few decades, see
for instance \cite{DengH1999,Panz1999}.  Recently,
event-triggered control for stochastic systems has also been drawing great attention, see, e.g., \cite{deng20191,deng20192,liuyg2019,liuyg2020}.
However, the mean power of white noise is infinite so that
SDE is not always appropriate for describing  many practical
systems, especially those in engineering \cite{randomsystem1}.
Random differential equations (RDEs) have been widely used as the dynamic
models, which are driven by second moment processes whose mean power is bounded, see, e.g., \cite{Khasminskii1980,Khalil1978}. It was not until 2014 that
the existence and uniqueness of solutions and stability
analysis of RDEs were developed in \cite{wuzjTAC2014}.
A series of relevant advancements occurred subsequently, see e.g., \cite{randomsystem1,YaoLAutomatica2020,zhanghTAC2021}. For example, the trajectory tracking
control for Lagrange systems driven by
second moment processes was addressed  in \cite{randomsystem1};
 The existence and uniqueness of the global solution, stability, and
 the adaptive output feedback regulation control design
 for random nonlinear time-delay systems were investigated in \cite{YaoLAutomatica2020};
 The event-triggered adaptive tracking control for RDE systems with coexisting parametric uncertainties and severe nonlinearities
was addressed in \cite{zhanghTAC2021}.
 More research progresses on random systems in the field
of control can be found in the recent review paper \cite{xiR2022}.
Recently, linear and nonlinear event-triggered ESOs have been designed for the open-loop of a class of uncertain
random nonlinear systems driven by bounded noise and colored noise which are representative second moment processes, and the mean square and almost sure convergence of proposed ESOs
have been presented with rigorous theoretical analysis \cite{zhwueventESO}.

A prominent element for the ETM-based observer and controller to be workable is to prevent the Zeno phenomenon (i.e., the triggering
conditions are satisfied infinite times in finite time). This leads to enormous difficulty
in developing event-triggered control for stochastic/random systems, because the execution/sampling times and the
inter-execution times are both random and then a positive constant lower
bound for the random inter-execution times is difficult or even impossible to be obtained.
In aforementioned literature like \cite{deng20191,deng20192,liuyg2019,liuyg2020}, novel ETM with dwell time (time-regularization)
or periodic ETM was designed, with the Zeno phenomenon be directly excluded but the
theoretical analysis be much more sophisticated.

 In this paper, we  develop the
 event-triggered ADRC controller
design and convergence analysis for a class of uncertain random nonlinear systems. The main contribution and novelty can be
summed up as follows: a) The controlled random nonlinear systems are subject to broad scale random
total disturbance involving the coupling of nonlinear unmodeled dynamics, external deterministic disturbance, bounded noise, and colored noise;
b) The event-triggered ADRC controller is designed with two feasible event-triggering mechanisms be proposed; c) Both the mean square and almost surely
practical convergence of the ADRC's closed-loop of the uncertain random nonlinear systems is presented with rigorous theoretical proofs.

This paper has the following structure. In Section \ref{Se2}, the problem is formulated, and some preliminaries are
introduced. In Section \ref{SetionIII}, the event-triggered ADRC controller is designed and the main results  are presented.
The theoretical proofs of main results are given in Section \ref{proofofmainresult}.
In Section \ref{sectionIV}, some numerical simulations are carried out to validate the theoretical results, followed up  with  concluding remarks in Section \ref{sectionV}.

\section{Problem formulation and preliminaries}\label{Se2}

Throughout the paper, the following notations are used.
$\mathbb{E}$ denotes the mathematical expectation;
$|Z|$ stands for the absolute value of a scalar $Z$, and
$\|Z\|$ denotes the 2-norm (or Euclidean norm) of a vector $Z$;
$\mathbb{I}_{n}$
denotes the $n$-dimensional identity matrix; $\bar{x}_{i}\triangleq(x_{1},\cdots,x_{i})$, $\hat{\bar{x}}_{i}\triangleq(\hat{x}_{1},\cdots,\hat{x}_{i})$,
$\hat{x}=\hat{\bar{x}}_{n+1}$;
 $\lambda_{\min}(Z)$  and  $\lambda_{\max}(Z)$ denote
the minimum eigenvalue and maximum  eigenvalue of a positive definite matrix $Z$;
$\mathbb{I}_{\Omega_{*}}$ denotes an indicator function
with the function value be $1$ in the domain $\Omega_{*}$ and be $0$ otherwise;
 $(\Omega,\mathcal{F},\mathbb{F}, P)$ represents
a complete filtered probability space with a filtration
$\mathbb{F}=\{\mathcal{F}_{t}\}_{t\geq 0}$, where
two mutually independent one-dimensional standard Brownian motions $B_{i}(t)\;(i=1,2)$ are defined.

In this paper, the event-triggered ADRC approach is addressed for the output-feedback stabilization  and
disturbance rejection for  a class of uncertain random nonlinear systems driven by bounded noise and colored noise as follows:
\begin{equation}\label{system1.2}
\left\{\begin{array}{l} \dot{x}_{1}(t)=x_{2}(t)+g_{1}(x_{1}(t)), \cr
\dot{x}_{2}(t)=x_{3}(t)+g_{2}(x_{1}(t),x_{2}(t)), \cr \hspace{1.2cm} \vdots \cr
\dot{x}_{n}(t)=f(t,x(t),w_{1}(t),w_{2}(t))+g_{n}(x(t))+u(t),\cr
y(t)=x_{1}(t),
\end{array}\right.
\end{equation}
where $x(t)=(x_{1}(t),\cdots,x_{n}(t))\in \mathbb{R}^{n}$,
$u(t)\in \mathbb{R}$, and $y(t)\in \mathbb{R}$ are the state,
control input, and  output measurement of system, respectively; $f:[0,\infty)\times \mathbb{R}^{n+2}\rightarrow \mathbb{R}$
is an unknown system function, and $g_{i}:\mathbb{R}^{i}\rightarrow \mathbb{R}\; (i=1,\cdots,n)$ are the known ones,
but $g_{i}(\bar{x}_{i}(t))\;(i=2,\cdots,n)$ still indicate unknown dynamics because of the unmeasurable state;
  $w_{1}(t)\triangleq\psi(t,B_{1}(t))$
 defined by an unknown function $\psi:[0,\infty)\times \mathbb{R}\rightarrow \mathbb{R}$ satisfying Assumption (A2) is the bounded noise,
and $w_{2}(t)$ is the colored noise that is the solution to the It\^{o}-type stochastic differential equation
(see, e.g., \cite[p.426]{colorednoise}, \cite[p.101]{mao}):
\begin{equation}\label{21equ}
dw_{2}(t)=-\rho_{1} w_{2}(t)dt+\rho_{1} \sqrt{2\rho_{2}}dB_{2}(t),
\end{equation}
where $\rho_{1}>0$ and $\rho_{2}>0$ are constants representing
the correlation time and the noise intensity, respectively,
 and the initial values $w_{2}(0)\in \mathbb{R}$.
 $\rho_{1}$, $\rho_{2}$ could be unknown but within known bounded intervals.

The significance of considering these two kinds of random noises can be expounded as follows.
It is known that white noise can be a characterization of random disturbances in practice,
which is generally understood as a stationary process with
zero mean and constant power spectral density.
However, white noise is not always workable in characterizing noise in many engineering applications
because the mean power of continuous-time white noise is infinite and  it has a correlation time of $0$ \cite{randomsystem1}.
The colored noise described by \dref{21equ} is a second moment process, which is with  uneven power spectral density
and bounded mean power. Compared to white noise, colored noise could be
more realistic in describing noise  when
the real processes have finite or even long correlation time and bounded mean power \cite{colorednoise,randomsystem1}.
More explanations and comparisons of white noise and colored noise and physical backgrounds of
colored noise can be found in \cite{Khalil1978,colorednoise,wuzjTAC2014,randomsystem1,xiR2022}, etc.
As for the bounded noise,  conventional deterministic disturbance is its special example by letting $w_{1}(t)=\psi(t)$
which is the function with regard to the time argument only, and frequently occurring bounded noises
like $\sin(t+B_{1}(t))$ and  $\cos(t+B_{1}(t))$ in practice (\cite{boundednoise1}) are the concerned
ones satisfying the following Assumption (A2).

It was shown in \cite{uniformlyobserverable}  that any uniform observable affine single-input single-output (SISO) nonlinear system can be transformed
into the lower triangular form \dref{system1.2} with $w_{i}(t)\equiv 0\;(i=1,2)$. Moreover, system \dref{system1.2}
covers the essential-integral-chain system with matched disturbance and uncertainty as a special case of $g_{i}(\cdot)=0\;(i=1,\cdots,n)$, which is the normalized form
to demonstrate the designing process and applications of ADRC.
Numerous engineering systems can be described by system \dref{system1.2}
 like aforementioned dynamics of DC-DC power converter \cite{DC-DC}, permanent magnet synchronous motor \cite{permanentmagnetsynchronousmotor}, and Delta robot \cite{robotcontorl}, etc.

 It should be pointed out that the convergence of
  ESO for the open-loop systems in \cite{zhwueventESO} doesn't directly indicates the convergence of ESO and stability
 of controlled systems in the ADRC's closed-loop. This is because in the proof of convergence of
 estimation error of the random total disturbance including system state for the open-loop systems, the state
 should be assumed to be bounded in a statistical sense beforehand, while the boundedness of
 the closed-loop state depends on the convergence of  ESO conversely.
The boundedness assumption of state for the open-loop systems can be regarded as a slowly varying condition
 besides the common structural one (like exact observability). This presupposition is unwanted when the ADRC controller is designed, which is not an easy theoretical task.
In addition, both the designs and theoretical analysis of this paper are largely distinguished from the deterministic counterpart in \cite{eventADRC2017,eventADRC2019,eventADRC2021}.

\section{Event-triggered ADRC design and main results}\label{SetionIII}
Define
\begin{eqnarray}
x_{n+1}(t)\triangleq f(t,x(t),w_{1}(t),w_{2}(t)),
\end{eqnarray}
which is a random total disturbance (extended state) involving the nonlinear coupling effects of the nonlinear unmodeled dynamics, external deterministic disturbance, bounded noise, and colored noise.
To estimate both unmeasurable state and random total disturbance, the event-triggered ESO
 is designed via the output and input of system (\ref{system1.2}) as follows:
 \begin{equation}\label{observer2.4}
\hspace{-0.1cm}\left\{\begin{array}{l}
\dot{\hat{x}}_{1}(t)=\hat{x}_{2}(t)+\lambda_{1}r\left(y(t_{k})-\hat{x}_{1}(t)\right)+g_{1}(\hat{x}_{1}(t)),
\crr
\dot{\hat{x}}_{2}(t)=\hat{x}_{3}(t)+\lambda_{2}r^{2}\left(y(t_{k})-\hat{x}_{1}(t)\right)+g_{2}(\hat{x}_{1}(t),\hat{x}_{2}(t)),\crr
\hspace{1.2cm} \vdots \crr
 \dot{\hat{x}}_{n}(t)=\hat{x}_{n+1}(t)+\lambda_{n}r^{n}\left(y(t_{k})-\hat{x}_{1}(t)\right)+g_{n}(\hat{\bar{x}}_{n}(t)) \cr\hspace{1.2cm}+u(t),  \crr
 \dot{\hat{x}}_{n+1}(t)=\lambda_{n+1}r^{n+1}\left(y(t_{k})-\hat{x}_{1}(t)\right),\; t\in [t_{k},t_{k+1}),
\end{array}\right.
\end{equation}
where $k\in\mathbb{Z}^{+}$, $t_{1}=0$, $\hat{x}_{i}(t)\;(i=1,\cdots,n+1)$ are the estimates of $x_{i}(t)$,
$r$ is the tuning gain, parameters $\lambda_{i}\;(i=1,\cdots,n+1)$ are selected to guarantee that the matrix
\begin{equation}\label{matricf2}\disp
H=
\begin{pmatrix}
-\lambda_{1}    & 1      &0      &  \cdots   &  0             \cr -\lambda_{2}
& 0 & 1   &  \cdots   &  0             \cr \cdots &\cdots & \cdots &
\cdots& \cdots  \cr -\lambda_{n}     & 0     & 0     &  \ddots   &  1 \cr
-\lambda_{n+1} & 0 &0 & \cdots & 0
\end{pmatrix}_{(n+1)\times (n+1)}
\end{equation}
is Hurwitz, and $t_{k}\;(k\in\mathbb{Z}^{+})$ are random execution times (stopping times) determined by the following event-triggering ETM
\begin{equation}\label{triggermechemi}
t_{k+1}=\inf\{t\geq t_{k}+\tau:\;|y(t)-y(t_{k})|\geq \kappa_{1} r^{-(n+\frac{1}{2})} \},
\end{equation}
with $\tau=\epsilon_{1} r^{-(2n+\frac{3}{2})}$ and $\epsilon_{1}$, $\kappa_{1}$ be free positive tuning parameters.
For any fixed $t\in[0,\infty)$, the last execution time of the ETM \dref{triggermechemi} before $t$ can be expressed as follows:
\begin{eqnarray}\label{fdfew}
\zeta_{t}=\max\{t_{k}:t_{k}\leq t, k\in \mathbb{Z}^{+}\}.
\end{eqnarray}
Therefore, $y(\zeta_{t})-y(t)$ can always represent $y(t_{k})-y(t),\;t\in [t_{k},t_{k+1})$ for any $k\in \mathbb{Z}^{+}$. Since $t_{k+1}-t_{k}\geq \tau$, the maximum
of execution times before $t$ is $[\frac{t}{\tau}]+1$ for almost every sample path. For $k=1,\cdots,[\frac{t}{\tau}]+1$, we define
\begin{eqnarray}
&& \Omega_{k}=\{\zeta_{t}=t_{k}\}, \; \Omega_{k,\tau}=\{\zeta_{t}=t_{k}\;\mbox{and}\; t\leq t_{k}+\tau\}.
\end{eqnarray}
Consequently, $\Omega$ can be expressed as the union of a set of mutually disjoint subsets as
$\Omega=\bigcup^{[\frac{t}{\tau}]+1}_{k=1} \Omega_{k}$ for any fixed $t\in[0,\infty)$.
Based on the event-triggered ESO \dref{observer2.4}, the event-triggered ADRC controller is
designed as
\begin{eqnarray}\label{triggeredadrc}
u(t)=\sum^{n}_{i=1}\theta^{n+1-i}c_{i}\hat{x}_{i}(t^{*}_{l})-\hat{x}_{n+1}(t^{*}_{l}), \; t\in [t^{*}_{l},t^{*}_{l+1}),
\end{eqnarray}
where $t^{*}_{1}=0$, $\theta\geq 1$ is to be specified later, $c_{i}\;(i=1,\cdots,n)$ are the control gains chosen such that the matrix
\begin{equation}\label{fe5gffd}
J=
\begin{pmatrix}
0     & 1      &0      &  \cdots   &  0             \cr 0     & 0 &
1     &  \cdots   &  0             \cr \cdots &\cdots & \cdots &
\cdots& \cdots  \cr 0     & 0     & 0     &  \ddots   &  1 \cr c_{1}
& c_{2} &\cdots & c_{n-1} & c_{n}
\end{pmatrix}_{n\times n}
\end{equation}
is Hurwitz, and $t^{*}_{l}\;(l\in\mathbb{Z}^{+})$ are random execution times (stopping times) determined by the following ETM
\begin{equation}\label{triggermechem2}
t^{*}_{l+1}=\inf\{t\geq t^{*}_{l}+\upsilon:\;\sum^{n+1}_{i=1}|\hat{x}_{i}(t)-\hat{x}_{i}(t^{*}_{l})|\geq \frac{\kappa_{2}}{r^{\frac{1}{2}}}\},
\end{equation}
with $l\in\mathbb{Z}^{+}$, $\upsilon=\epsilon_{2}r^{-(\frac{2n}{3}+\frac{5}{3})}$, and $\epsilon_{2}$, $\kappa_{2}$ be any free positive tuning parameters.
The  event-triggered ADRC controller \dref{triggeredadrc} is comprised of an output-feedback controller
$\sum^{n}_{i=1}\theta^{n+1-i}c_{i}\hat{x}_{i}(t^{*}_{l})$  and a disturbance rejection component $-\hat{x}_{n+1}(t^{*}_{l})$
under the ETM \dref{triggermechem2}, which rejects the random total disturbance in an active way
but not the passive one.

Similarly, for any fixed $t\in[0,\infty)$, the last execution time of the ETM \dref{triggermechem2} before $t$ is
\begin{eqnarray}\label{wtt}
 \varpi_{t}\triangleq\max\{t^{*}_{l}:t^{*}_{l}\leq t, l\in \mathbb{Z}^{+}\}.
\end{eqnarray}
Since $t^{*}_{l+1}-t^{*}_{l}\geq \upsilon $, the maximum
of execution times before $t$ is $[\frac{t}{\upsilon}]+1$ for almost every sample path. For $l=1,\cdots,[\frac{t}{\upsilon}]+1$, we define
\begin{eqnarray}
&& \Omega^{*}_{l}=\{\varpi_{t}=t^{*}_{l}\}, \; \Omega^{*}_{l,\upsilon}=\{\varpi_{t}=t^{*}_{l}\;\mbox{and}\; t\leq t^{*}_{l}+\upsilon\}.
\end{eqnarray}
Then $\Omega$ can also be represented as the union of a set of mutually disjoint subsets as
$\Omega=\bigcup^{[\frac{t}{\upsilon}]+1}_{l=1} \Omega^{*}_{l}$ for any fixed $t\in[0,\infty)$.
In addition, we also set
\begin{eqnarray}
\Omega_{\tau}=\bigcup^{[\frac{t}{\tau}]+1}_{k=1} \Omega_{k,\tau},\;\;\; \Omega^{*}_{\upsilon}=\bigcup^{[\frac{t}{\upsilon}]+1}_{l=1} \Omega^{*}_{l,\upsilon}.
\end{eqnarray}

\begin{remark}
With regard to the ETM \dref{triggermechemi} (or \dref{triggermechem2}),
each inter-execution
interval $[t_{k},t_{k+1})$ (or $[t^{*}_{l},t^{*}_{l+1})$) is of two-stage: once the execution emerges, the trigger
will still cease in $[t_{k},t_{k}+\tau)$ (or $[t^{*}_{l},t^{*}_{l}+\upsilon)$); and the event triggering condition is continuously evaluated
after the time instant $t_{k}+\tau$ (or $t^{*}_{l}+\upsilon$) to determine the next execution time $t_{k+1}$ (or $t^{*}_{l+1}$).
This means that
 each inter-execution time  is clearly not less than $\tau$ (or $\upsilon$),
so that the Zeno phenomenon can be naturally avoided. However, the convergence analysis of ADRC's closed-loop
under these ETMs would be more complex.
\end{remark}

To guarantee the mean square and almost surely practical convergence  of the resulting closed-loop of system \dref{system1.2} under the  event-triggered ADRC controller \dref{triggeredadrc},
the following assumptions are required.

{\bf  Assumption (A1).}
  $f(\cdot)$  has first-order and second-order continuous partial derivatives  with regard to the arguments $(t,x)$ and
$(w_{1},w_{2})$, respectively, and there are known constants  $\alpha_{j}>0
\;(j=1,2,3,4)$, $L_{i}\geq 0$, and  functions $\varphi_{j}\in
C(\mathbb{R};\mathbb{R}^{+})\;(j=1,2)$ such that for all $t\ge 0$, $x \in \mathbb{R}^{n}$,
$w_{1}\in \mathbb{R}$, $w_{2}\in \mathbb{R}$, $i=1,\cdots,n$, there holds
\begin{eqnarray*}\label{2.40}
\begin{array}{ll}
\left|f(t,x,w_{1},w_{2})\right|+\left|\frac{\partial f(t,x,w_{1},w_{2})}{\partial t}\right|
\cr\leq \alpha_{1}+\alpha_{2}\|x\|+\alpha_{3}|w_{2}|+\varphi_{1}(w_{1}),
\end{array}
\end{eqnarray*}
\begin{eqnarray*}
\begin{array}{ll}
 \disp\sum^{n}_{i=1}\left|\frac{\partial f(t,x,w_{1},w_{2})}{\partial x_{i}}\right|+
\left|\frac{\partial f(t,x,w_{1},w_{2})}{\partial w_{1}}\right|\cr \disp
+ \left|\frac{\partial^{2} f(t,x,w_{1},w_{2})}{\partial w^{2}_{1}}\right|+\left|\frac{\partial f(t,x,w_{1},w_{2})}{\partial w_{2}}\right|
\cr \disp+ \left|\frac{\partial^{2} f(t,x,w_{1},w_{2})}{\partial w^{2}_{2}}\right|\leq \alpha_{4}+\varphi_{2}(w_{1}),
\end{array}
\end{eqnarray*}
\begin{eqnarray*}
\left|g_{i}(\bar{x}_{i})-g_{i}(\hat{\bar{x}}_{i})\right|\leq L_{i}\|\bar{x}_{i}-\hat{\bar{x}}_{i}\|,\; g_{i}(\underbrace{0,\cdots,0}_{i})=0.
\end{eqnarray*}

{\bf Assumption (A2).} The function $\psi(t,\vartheta):[0,\infty)\times \mathbb{R} \rightarrow \mathbb{R}$
has first order and second order continuous partial derivatives with regard to
the arguments $t$ and $\vartheta$, respectively, and there exists a known constant $\alpha_{5}>0$, such that for all
  $t\geq 0$, $\vartheta\in \mathbb{R}$,
\begin{eqnarray*}\label{assumptiodn}
\left|\psi(t,\vartheta)\right|+\left|\frac{\partial \psi(t,\vartheta)}{\partial t}\right|
+\left|\frac{\partial \psi(t,\vartheta)}{\partial \vartheta}\right|+
\frac{1}{2}\left|\frac{\partial^{2} \psi(t,\vartheta)}{\partial \vartheta^{2}}\right|\leq \alpha_{5}.
\end{eqnarray*}
\begin{remark}
Because the random total disturbance is to be estimated by the event-triggered ESO and compensated in the event-triggered ADRC's closed-loop,
 its rate of variation (stochastic differential in \dref{stocahsticdifferential}) is naturally be required to be bounded or
linear growth with respect to the closed-loop states, which is guaranteed by  Assumptions (A1)-(A2).
\end{remark}
In the following main results and their proofs, the following symbols are used:
\begin{eqnarray}\label{zdefintion}
&& z_{i}(t)=r^{n+1-i}(x_{i}(t)-\hat{x}_{i}(t)),\;i=1,\cdots,n+1,\cr&&
z(t)=(z_{1}(t),\cdots,z_{n+1}(t)),\; c_{n+1}=1.
\end{eqnarray}

Set
\begin{eqnarray}\label{414df}
&&\hspace{-0.1cm}
\Lambda_{1}(\upsilon)= \frac{10(n+1)\theta^{2n}(\disp\max_{1\leq i\leq n+1}c^{2}_{i}+4\alpha^{2}_{2})\upsilon^{2}}{(1-\upsilon\theta^{n}\disp\max_{1\leq i\leq n+1}|c_{i}|)^{2}}, \cr &&\hspace{-0.1cm}
 \Lambda_{2}(\upsilon)=\frac{10[n(1+\disp\sum^{n}_{i=1}L^{2}_{i})+4n\alpha^{2}_{2}]\upsilon}{(1-\upsilon\theta^{n}\disp\max_{1\leq i\leq n+1}|c_{i}|)^{2}}, \cr&&\hspace{-0.1cm}
\Lambda_{3}(\upsilon,\tau,r)=\frac{10(n+1)(1+L^{2}_{1})\upsilon^{2}(\upsilon+\tau)r^{2(n+1)}\disp\sum^{n+1}_{i=1}\lambda^{2}_{i}}{(1-\upsilon\theta^{n}\disp\max_{1\leq i\leq n+1}|c_{i}|)^{2}}, \cr&&\hspace{-0.1cm}
\Lambda_{4}(\upsilon)=\frac{10(n+1)\disp\theta^{2n}\max_{1\leq i\leq n+1}c^{2}_{i}\upsilon^{2}}{(1-\upsilon\theta^{n}\disp\max_{1\leq i\leq n+1}|c_{i}|)^{2}}, \cr&&\hspace{-0.1cm}
\Lambda_{5}(\upsilon,r)=\frac{10[n(1+\disp\sum^{n}_{i=1}\frac{L^{2}_{i}}{r^{2(n+1-i)}})+(n+1)r^{2}\sum^{n+1}_{i=1}\lambda^{2}_{i}]\upsilon }{(1-\upsilon\theta^{n}\disp\max_{1\leq i\leq n+1}|c_{i}|)^{2}}, \cr&&\hspace{-0.1cm}
\Lambda_{6}(\upsilon,r)=\frac{10\theta^{2n}\upsilon^{2}}{(1-\upsilon\theta^{n}\disp\max_{1\leq i\leq n+1}|c_{i}|)^{2}}\bigg\{(4+8n)[\alpha^{2}_{1}+\cr&&\hspace{-0.1cm}\alpha^{2}_{3}\sup_{t\in [-\upsilon,\infty)}\mathbb{E}|w_{2}(t)|^{2}
+\sup_{t\in [-\upsilon,\infty)}(\varphi_{1}(w_{1}(t)))^{2}]\cr&&\hspace{-0.1cm}+r(n+1)\kappa^{2}_{1}\sum^{n+1}_{i=1}\lambda^{2}_{i}\bigg\}, \cr&&\hspace{-0.1cm}
x(t)\triangleq x(0),\; t\in [-\upsilon-\tau,0],\cr &&\hspace{-0.1cm} z(t)\triangleq z(0), w_{i}(t)\triangleq w_{i}(0), t\in [-\upsilon,0],i=1,2,
\end{eqnarray}
where tuning parameters $\upsilon$, $\theta$, $c_{i}$ are chosen such that $\upsilon\theta^{n}\max_{1\leq i\leq n+1}\{|c_{i}|\}<1$.

The sampling error
$\hat{x}_{i}(t)-\hat{x}_{i}(\varpi_{t})$ on $\Omega^{*}_{v}$ is provided in Lemma \ref{errorlemma}, necessary for the convergence
analysis of the closed-loop systems later.
\begin{lemma}\label{errorlemma}
Suppose that Assumptions (A1)-(A2) hold, and  tuning parameters $\upsilon$, $\theta$, $c_{i}$ are chosen such that $\upsilon\theta^{n}\max_{1\leq i\leq n+1}\{|c_{i}|\}<1$, then for all $t\geq 0$, there holds
\begin{eqnarray}
&&\hspace{-0.1cm}\mathbb{E}[\sum^{n+1}_{i=1}|\hat{x}_{i}(t)-\hat{x}_{i}(\varpi_{t})|\mathbb{I}_{\Omega^{*}_{\upsilon}}]^{2}\cr &&\hspace{-0.1cm}\leq  \Lambda_{1}(\upsilon)\mathbb{E}[\|x(t)\|^{2}\mathbb{I}_{\Omega^{*}_{\upsilon}}]+\Lambda_{2}(\upsilon)\int^{t}_{t-\upsilon}\mathbb{E}[\|x(s)\|^{2}\mathbb{I}_{\Omega^{*}_{\upsilon}}]ds\cr&&\hspace{-0.1cm}
+\Lambda_{3}(\upsilon,\tau,r)\int^{t}_{t-\upsilon-\tau}\mathbb{E}[\|x(s)\|^{2}\mathbb{I}_{\Omega^{*}_{\upsilon}}]ds+\Lambda_{4}(\upsilon)\mathbb{E}[\|z(t)\|^{2}\mathbb{I}_{\Omega^{*}_{\upsilon}}]\cr&&\hspace{-0.1cm}
+\Lambda_{5}(\upsilon,r)\int^{t}_{t-\upsilon}\mathbb{E}[\|z(s)\|^{2}\mathbb{I}_{\Omega^{*}_{\upsilon}}]ds+
\Lambda_{6}(\upsilon,r),
\end{eqnarray}
where $\Lambda_{i}\;(i=1,2,\cdots,6)$ are specified in \dref{414df}.

\end{lemma}
\textbf{Proof.} See ``Proof of Lemma \ref{errorlemma}" in Section \ref{proofofmainresult}.

Let $Q_{1}\in \mathbb{R}^{n\times n}$ and $Q_{2}\in \mathbb{R}^{(n+1)\times (n+1)}$
be   the respective unique positive definite matrices solutions of
the Lyapunov equations
\begin{eqnarray}\label{Q1E}
Q_{1}J+J^{\top}Q_{1}=-\mathbb{I}_{n}
\end{eqnarray}
and
\begin{eqnarray}\label{Q2E}
Q_{2}H+H^{\top}Q_{2}=-\mathbb{I}_{n+1},
\end{eqnarray}
where $J$ and $H$ are given in \dref{fe5gffd} and \dref{matricf2}, respectively.

The mean square practical convergence and almost surely practical one of
the closed-loop of system \dref{system1.2} under the event-triggered ADRC controller \dref{triggeredadrc} is summarized up
as the following Theorem \ref{theorem34.1}.
\begin{theorem}\label{theorem34.1}
Suppose that Assumptions (A1)-(A2) hold, $\theta$ is chosen such that $\theta>2\lambda_{\max}(Q_{1})\disp\sum^{n}_{i=1}L_{i}$ and $\upsilon\theta^{n}\disp\max_{1\leq i\leq n+1}\{|c_{i}|\}<1$.
Then, there exists a known constant $r_{*}\geq 1$, such that for any initial values $x(0)\in \mathbb{R}^{n}$, $\hat{x}(0)\in \mathbb{R}^{n+1}$, $w_{2}(0)\in \mathbb{R}$,
and any $r\geq r_{*}$,  there exists a unique global solution $(x(t),\hat{x}(t))$ to the closed-loop of system \dref{system1.2} under the  event-triggered ADRC controller \dref{triggeredadrc} that satisfies
\begin{eqnarray}\label{estimationerror}
&&\mbox{(i)}\;\;\;\mathbb{E}|x_{i}(t)-\hat{x}_{i}(t)|^2\leq \frac{M}{r^{2n+3-2i}},\; i=1,\cdots,n+1,\cr && \disp
\mbox{(ii)}\;\;|x_{i}(t)-\hat{x}_{i}(t)| \leq \frac{M_{\omega}}{r^{\frac{2n+3-2i}{2}}} \; \mbox{a.s.} \; i=1,\cdots,n+1,\cr &&\disp
\mbox{(iii)}\;\sum^{n}_{i=1}\mathbb{E}|x_{i}(t)|^2\leq \frac{M}{r}, \cr &&\disp
\mbox{(iv)}\;\; |x_{i}(t)|\leq \frac{M_{\omega}}{r^{\frac{1}{2}}} \; \mbox{a.s.} \; i=1,\cdots,n,
\end{eqnarray}
for all $t\geq t_{r}$ with $t_{r}$ be an $r$-dependent constant, where $M>0$ and $M_{\omega}>0$ are a constant and a random variable  independent of $r$, respectively.
\end{theorem}

\textbf{Proof.} See ``Proof of Theorem \ref{theorem34.1}" in Section \ref{proofofmainresult}.


\begin{remark}
Compared to the convergence of linear event-triggered ESO in \cite{zhwueventESO}, both the mean square practical convergence and almost surely practical one of the
event-triggered ESO in the closed-loop are obtained without requiring the boundedness of state as a presupposition.
The concerning convergence in aforementioned main results mean that the estimation error and the bound of the closed-loop state in the mean square and almost sure sense
can be as small as we want in $[t_{r},\infty)$, provided that the gain $r$ is tuned to be large enough, which requires the $r$-dependent event-triggered ESO-based controller
to be designed accordingly.
\end{remark}


\section{Proofs of the main results}\label{proofofmainresult}

The proofs of the main results are given in this section.

{\bf Proof of Lemma \ref{errorlemma}.} The following procedures in the proof of Lemma \ref{errorlemma} are implemented
for any fixed $t\in [0,\infty)$. Let
\begin{eqnarray}\label{df32}
\delta_{i}(t)=\hat{x}_{i}(t)-\hat{x}_{i}(\varpi_{t}),
\end{eqnarray}
where $\varpi_{t}$ is defined as that in \dref{wtt}. $c_{n+1}$, $z_{i}(t)$, and $z(t)$ are
specified in \dref{zdefintion}.
For $\omega\in\Omega_{k,\tau}$, we have $\zeta_{t}\geq t-\tau$. Thus, by the ETM \dref{triggermechemi}, for all $t\geq 0$,
\begin{eqnarray*}
 &&|y(\zeta_{t})-y(t)|\cr&&=\sum^{[\frac{t}{\tau}]+1}_{k=1}|y(\zeta_{t})-y(t)|\mathbb{I}_{\Omega_{k,\tau}}+\sum^{[\frac{t}{\tau}]+1}_{k=1}|y(\zeta_{t})-y(t)|\mathbb{I}_{\Omega_{k}\setminus\Omega_{k,\tau}}
 \cr&& \leq \int^{t}_{t-\tau}|x_{2}(s)|ds+L_{1}\int^{t}_{t-\tau}|x_{1}(s)|ds+\kappa_{1}r^{-(n+\frac{1}{2})}.
\end{eqnarray*}
Therefore, for $i=1,\cdots,n+1$, it follows that
\begin{eqnarray*}
 &&\int^{t}_{t-\upsilon}|\lambda_{i}r^{i}(y(\zeta_{s})-y(s))|ds\cr&&=
 |\lambda_{i}|\upsilon r^{i}\int^{t}_{t-\upsilon-\tau}|x_{2}(s)|ds+ |\lambda_{i}|L_{1}\upsilon r^{i}\int^{t}_{t-\upsilon-\tau}|x_{1}(s)|ds\cr&&+|\lambda_{i}|\kappa_{1}\upsilon r^{-(n+\frac{1}{2}-i)}.
\end{eqnarray*}
For $\omega\in\Omega^{*}_{l,v}$, we have $\varpi_{t}\geq t-\upsilon$.
Then, for any $\omega\in\Omega^{*}_{l,v}$ and $i=1,\cdots,n-1$, it is obtained that
\begin{eqnarray*}
&& |\delta_{i}(t)|\cr&&\leq
\int^{t}_{t-\upsilon}|\hat{x}_{i+1}(s)+\lambda_{i}r^{i}(y(\zeta_{s})-\hat{x}_{1}(s))+g_{i}(\hat{\bar{x}}_{i}(s))|ds \cr&&\leq
\int^{t}_{t-\upsilon}\frac{|z_{i+1}(s)|}{r^{n-i}}ds+\int^{t}_{t-\upsilon}|x_{i+1}(s)|ds+
|\lambda_{i}|\upsilon r^{i}\int^{t}_{t-\upsilon-\tau}\cr&&|x_{2}(s)|ds+ |\lambda_{i}|L_{1}\upsilon r^{i}\int^{t}_{t-\upsilon-\tau}|x_{1}(s)|ds+|\lambda_{i}|\kappa_{1}\upsilon r^{-(n+\frac{1}{2}-i)}\cr&&+
\int^{t}_{t-\upsilon}\frac{|\lambda_{i}z_{1}(s)|}{r^{n-i}}ds
+\frac{L_{i}}{r^{n+1-i}}\int^{t}_{t-\upsilon}\|z(s)\|ds\cr&&+L_{i}\int^{t}_{t-\upsilon}\|x(s)\|ds,
\end{eqnarray*}
\begin{eqnarray*}
&& |\delta_{n}(t)|\cr&& \leq \sum^{n+1}_{i=1}\upsilon\theta^{n+1-i}|c_{i}\hat{x}_{i}(\varpi_{t})|
+\int^{t}_{t-\upsilon}|\lambda_{n}r^{n}(y(\zeta_{s})-\hat{x}_{1}(s))|ds\cr&&  +\int^{t}_{t-\upsilon}|\hat{x}_{n+1}(s)|ds+\int^{t}_{t-\upsilon}|g_{n}(\hat{\bar{x}}_{n}(t))|ds
\cr&& \leq \sum^{n+1}_{i=1}\upsilon \theta^{n+1-i}|c_{i}\delta_{i}(t)|+\sum^{n+1}_{i=1}\upsilon\theta^{n+1-i}\frac{|c_{i}z_{i}(t)|}{r^{n+1-i}}\cr&& +\sum^{n+1}_{i=1}\upsilon\theta^{n+1-i}|c_{i}x_{i}(t)|
+|\lambda_{n}|\upsilon r^{n}\int^{t}_{t-\upsilon-\tau}|x_{2}(s)|ds\cr&& +|\lambda_{n}|L_{1}\upsilon r^{n}\int^{t}_{t-\upsilon-\tau}|x_{1}(s)|ds+|\lambda_{n}|\kappa_{1}\upsilon r^{-\frac{1}{2}}\cr&& +\int^{t}_{t-\upsilon}|\lambda_{n}z_{1}(s)|ds
+\int^{t}_{t-\upsilon}|z_{n+1}(s)|ds+\int^{t}_{t-\upsilon}|x_{n+1}(s)|ds\cr&& +\frac{L_{n}}{r}\int^{t}_{t-\upsilon}\|z(s)\|ds+L_{n}\int^{t}_{t-\upsilon}\|x(s)\|ds,
\end{eqnarray*}
and
\begin{eqnarray*}
&&|\delta_{n+1}(t)|
\leq |\lambda_{n+1}|\upsilon r^{n+1}\int^{t}_{t-\upsilon-\tau}|x_{2}(s)|ds\cr&&+ |\lambda_{n+1}|L_{1}\upsilon r^{n+1}\int^{t}_{t-\upsilon-\tau}|x_{1}(s)|ds +|\lambda_{n+1}|\kappa_{1}\upsilon r^{\frac{1}{2}}\cr&&+r\int^{t}_{t-\upsilon}|\lambda_{n+1}z_{1}(s)|ds.
\end{eqnarray*}
Therefore, for any $\omega\in\Omega^{*}_{l,v}$, we have
\begin{eqnarray*}
&&\sum^{n+1}_{i=1}|\delta_{i}(t)|\cr&&\leq\frac{1}{1-\upsilon\theta^{n}\disp\max_{1\leq i\leq n+1}|c_{i}|}\bigg\{\sum^{n}_{i=1}\int^{t}_{t-\upsilon}\frac{|z_{i+1}(s)|}{r^{n-i}}ds
\cr&&+\sum^{n}_{i=1}\int^{t}_{t-\upsilon}|x_{i+1}(s)|ds+\sum^{n+1}_{i=1}|\lambda_{i}|\upsilon r^{i}\int^{t}_{t-\upsilon-\tau}|x_{2}(s)|ds\cr&&+
\sum^{n+1}_{i=1} |\lambda_{i}|L_{1}\upsilon r^{i}\int^{t}_{t-\upsilon-\tau}|x_{1}(s)|ds+\sum^{n+1}_{i=1}\int^{t}_{t-\upsilon}\frac{|\lambda_{i}z_{1}(s)|}{r^{n-i}}ds\cr&&+\sum^{n+1}_{i=1}\upsilon\theta^{n+1-i}\frac{|c_{i}z_{i}(t)|}{r^{n+1-i}}+\sum^{n+1}_{i=1}\upsilon\theta^{n+1-i}|c_{i}x_{i}(t)|
\cr&&+\sum^{n+1}_{i=1}|\lambda_{i}|\kappa_{1}\upsilon r^{-(n+\frac{1}{2}-i)}+\sum^{n}_{i=1}\frac{L_{i}}{r^{n+1-i}}\int^{t}_{t-\upsilon}\|z(s)\|ds\cr&&+\sum^{n}_{i=1}L_{i}\int^{t}_{t-\upsilon}\|x(s)\|ds\bigg\}.
\end{eqnarray*}
These, together with Assumptions (A1)-(A2) and the inequality $\disp(\sum^{m}_{i=1}a_{i})^{2}\leq m \sum^{m}_{i=1}a^{2}_{i}$ for $a_{i}\geq 0$ and $m\in \mathbb{Z}^{+}$, further yield that
\begin{eqnarray}
&& \hspace{-0.2cm}\mathbb{E}[\sum^{n+1}_{i=1}|\delta_{i}(t)|\cdot\mathbb{I}_{\Omega^{*}_{v}}]^{2}
\cr&&\hspace{-0.2cm}\leq \frac{10}{(1-\upsilon\theta^{n}\disp\max_{1\leq i\leq n+1}|c_{i}|)^{2}}\sum^{[\frac{t}{\upsilon}]+1}_{l=1}\mathbb{E}\bigg\{\bigg\{n\upsilon(1+\sum^{n}_{i=1}L^{2}_{i})\cdot\cr&&\hspace{-0.2cm}\int^{t}_{t-\upsilon}\|x(s)\|^{2}ds+n\upsilon(1+\sum^{n}_{i=1}\frac{L^{2}_{i}}{r^{2(n+1-i)}})\int^{t}_{t-\upsilon}\|z(s)\|^{2}
ds\cr&&\hspace{-0.2cm}
+(n+1)(1+L^{2}_{1})\upsilon^{2}(\upsilon+\tau)r^{2(n+1)}\sum^{n+1}_{i=1}\lambda^{2}_{i}\int^{t}_{t-\upsilon-\tau}\|x(s)\|^{2}\cr&&\hspace{-0.2cm}ds
+(n+1)\upsilon r^{2}\sum^{n+1}_{i=1}\lambda^{2}_{i}\int^{t}_{t-\upsilon}\|z(s)\|^{2}ds
+(n+1)\upsilon^{2}\theta^{2n}\cdot\cr&&\hspace{-0.2cm}\max_{1\leq i\leq n+1}c^{2}_{i}(\|x(t)\|^{2}+\|z(t)\|^{2})
+4n\upsilon\int^{t}_{t-\upsilon}[\alpha^{2}_{1}+\alpha^{2}_{2}\|x(s)\|^{2}\cr&&\hspace{-0.2cm}+\alpha^{2}_{3}|w_{2}(s)|^{2}+(\varphi_{1}(w_{1}(s)))^{2}] ds+4(n+1)\upsilon^{2}\theta^{2n} [\alpha^{2}_{1}\cr&&\hspace{-0.2cm}+\alpha^{2}_{2}\|x(t)\|^{2}+\alpha^{2}_{3}|w_{2}(t)|^{2}+(\varphi_{1}(w_{1}(t)))^{2}]\cr&&\hspace{-0.2cm}+\upsilon^{2}r(n+1)\kappa^{2}_{1}\sum^{n+1}_{i=1}\lambda^{2}_{i}\bigg\}\mathbb{I}_{\Omega^{*}_{l,\upsilon}}
\bigg\} \cr &&\hspace{-0.2cm}\leq  \Lambda_{1}(\upsilon)\mathbb{E}[\|x(t)\|^{2}\mathbb{I}_{\Omega^{*}_{v}}]+\Lambda_{2}(\upsilon)\int^{t}_{t-\upsilon}\mathbb{E}[\|x(s)\|^{2}\mathbb{I}_{\Omega^{*}_{v}}]ds\cr&&\hspace{-0.2cm}
+\Lambda_{3}(\upsilon,\tau,r)\int^{t}_{t-\upsilon-\tau}\mathbb{E}[\|x(s)\|^{2}\mathbb{I}_{\Omega^{*}_{v}}]ds+\Lambda_{4}(\upsilon)\mathbb{E}[\|z(t)\|^{2}\cdot\mathbb{I}_{\Omega^{*}_{v}}]\cr&&\hspace{-0.2cm}
+\Lambda_{5}(\upsilon,r)\int^{t}_{t-\upsilon}\mathbb{E}[\|z(s)\|^{2}\mathbb{I}_{\Omega^{*}_{\upsilon}}]ds+
\Lambda_{6}(\upsilon,r),
\end{eqnarray}
where $\Lambda_{i}\;(i=1,\cdots,6)$ are specified in \dref{414df}.
This ends the proof.
\hfill $\Box$

{\bf Proof of Theorem \ref{theorem34.1}.}
Set
\begin{equation}
\left\{\begin{array}{l}
\varrho_{i}(t)=\theta^{n-i}x_{i}(t), \;\;
\varrho(t)=(\varrho_{1}(t),\cdots,\varrho_{n}(t)), \cr
\Xi_{i}(t)=g_{i}(\bar{x}_{i}(t))-g_{i}(\hat{\bar{x}}_{i}(t)),\;i=1,\cdots,n,\cr
\varsigma(t)=r^{n}(y(\zeta_{t})-y(t)),\cr
\delta(t)=(\delta_{1}(t),\cdots,\delta_{n+1}(t)),
\end{array}\right.
\end{equation}
where $\delta_{i}(t)\;(i=1,\cdots,n+1)$ are defined as those in \dref{df32}. $c_{n+1}$, $z_{i}(t)$, and $z(t)$ are
specified in \dref{zdefintion}.
Applying It\^{o}'s formula to the random total disturbance $f(t,x(t),w_{1}(t),w_{2}(t))$ with respect to
$t$ along system \dref{system1.2} and \dref{21equ} with the event-triggered ADRC controller designed in \dref{triggeredadrc}, it is obtained that
\begin{eqnarray}\label{stocahsticdifferential}
&&\hspace{-0.2cm} dx_{n+1}(t)=\bigg\{\frac{\partial f(t,x(t),w_{1}(t),w_{2}(t))}{\partial
t}\cr &&\hspace{-0.2cm}+ \sum\limits^{n-1}_{i=1}\frac{\partial
f(t,x(t),w_{1}(t),w_{2}(t))}{\partial x_{i}}[x_{i+1}(t)+g_{i}(\bar{x}_{i}(t))]
\cr &&\hspace{-0.2cm} +\frac{\partial f(t,x(t),w_{1}(t),w_{2}(t))}{\partial x_{n}}\cdot[x_{n+1}(t)+g_{n}(x(t))-\hat{x}_{n+1}(\varpi_{t})\cr &&\hspace{-0.2cm}+\sum^{n}_{i=1}\theta^{n+1-i}c_{i}\hat{x}_{i}(\varpi_{t})]dt +\frac{\partial
f(t,x(t),w_{1}(t),w_{2}(t))}{\partial
w_{1}}\cdot\cr&&\hspace{-0.2cm}[\frac{\partial \psi(t,B_{1}(t))}{\partial t}+\frac{1}{2}\frac{\partial^{2} \psi(t,B_{1}(t))}{\partial \vartheta^{2}}]\cr &&\hspace{-0.2cm} +
\frac{1}{2}\frac{\partial^{2}
f(t,x(t),w_{1}(t),w_{2}(t))}{\partial w^{2}_{1}}(\frac{\partial \psi(t,B_{1}(t))}{\partial \vartheta})^{2}\cr&&\hspace{-0.2cm}-
 \frac{\partial f(t,x(t),w_{1}(t),w_{2}(t))}{\partial w_{2}}\rho_{1}w_{2}(t)\cr&&\hspace{-0.2cm}+\frac{\partial^{2} f(t,x(t),w_{1}(t),w_{2}(t))}{\partial w^{2}_{2}}\rho^{2}_{1}\rho_{2}\bigg\}dt \cr&&\hspace{-0.2cm}+
 \frac{\partial f(t,x(t),w_{1}(t),w_{2}(t))}{\partial w_{1}}\frac{\partial \psi(t,B_{1}(t))}{\partial \vartheta}dB_{1}(t) \cr&&\hspace{-0.2cm}+
 \frac{\partial f(t,x(t),w_{1}(t),w_{2}(t))}{\partial w_{2}}\rho_{1}\sqrt{2\rho_{2}}dB_{2}(t)\cr&&\hspace{-0.2cm}\triangleq\Gamma_{1}(t)dt+\Gamma_{2}(t)dB_{1}(t)+\Gamma_{3}(t)dB_{2}(t),
\end{eqnarray}
where $\vartheta$ represents the second argument of the function $\psi(\cdot,\cdot)$, and $\varpi_{t}$ is defined in \dref{wtt}.
By Assumptions (A1)-(A2), there are known constants $\beta_{i}\;(i=1,\cdots,5)$ such that for all $t\geq 0$,
\begin{eqnarray}\label{55cite}
&&\mathbb{E}|\Gamma_{1}(t)|^{2}\leq \beta_{1}+\beta_{2}\mathbb{E}\|\varrho(t)\|^{2}+\beta_{3}\mathbb{E}\|z(t)\|^{2}+\beta_{4}\mathbb{E}\|\delta(t)\|^{2},\cr&&
|\Gamma_{2}(t)|^{2}+|\Gamma_{3}(t)|^{2}\leq \beta_{5}.
\end{eqnarray}
By Assumption (A1) about the function $g_{i}(\cdot)$, it follows that
\begin{eqnarray}\label{56cite}
|\Xi_{i}(t)|\leq \frac{L_{i}}{r^{n+1-i}}\|z(t)\|,\; \forall t\geq 0.
\end{eqnarray}
The closed-loop of system \dref{system1.2} under aforementioned event-triggered ADRC controller \dref{triggeredadrc} is equivalent to
\begin{equation}\label{closed-loopsystem}
\left\{\begin{array}{l}
d\varrho_{1}(t)=\theta\varrho_{2}(t)dt+\theta^{n-1}g_{1}(\frac{\varrho_{1}(t)}{\theta^{n-1}})dt,
\crr d\varrho_{2}(t)=\theta\varrho_{3}(t)dt+\theta^{n-2}g_{2}(\frac{\varrho_{1}(t)}{\theta^{n-1}},\frac{\varrho_{2}(t)}{\theta^{n-2}})dt,
\cr
\hspace{1.2cm}  \vdots \cr d\varrho_{n}(t)=\disp[\theta\sum^{n}_{i=1}c_{i}\varrho_{i}(t)-\sum^{n+1}_{i=1}\frac{\theta^{n+1-i}c_{i}z_{i}(t)}{r^{n+1-i}}\cr
-\disp\sum^{n+1}_{i=1}\theta^{n+1-i}c_{i}\delta_{i}(t)+g_{n}(\frac{\varrho_{1}(t)}{\theta^{n-1}},\cdots,\varrho_{n}(t))]dt,\crr
dz_{1}(t)=r[z_{2}(t)-\lambda_{1}z_{1}(t)]dt-\lambda_{1}r\varsigma(t)dt+r^{n}\Xi_{1}(t)dt,\cr
dz_{2}(t)=r[z_{3}(t)-\lambda_{2}z_{1}(t)]dt-\lambda_{2}r\varsigma(t)dt\cr\hspace{1.4cm}+r^{n-1}\Xi_{2}(t)dt, \cr\hspace{1.2cm}\vdots\cr
dz_{n}(t)=r[z_{n+1}(t)-\lambda_{n}z_{1}(t)]dt-\lambda_{n}r\varsigma(t)dt\cr\hspace{1.4cm}+r\Xi_{n}(t)dt,
 \cr
dz_{n+1}(t)=-r\lambda_{n+1}z_{1}(t)dt-\lambda_{n+1}r\varsigma(t)dt+\Gamma_{1}(t)dt\cr \hspace{1.8cm}+\Gamma_{2}(t)dB_{1}(t)+\Gamma_{3}(t)dB_{2}(t).
\end{array}\right.
\end{equation}
It follows from  (\ref{21equ}) that the colored noise $w_{2}(t)$ can be regarded as the
extended  state variable of (\ref{closed-loopsystem}), and the
bounded noise $w_{1}(t)$ is with deterministic bound satisfying Assumption (A2).
These,  together with Assumption (A1),  yield that the
local Lipschitz and the linear growth conditions are satisfied by
the drift term and diffusion one of (\ref{closed-loopsystem}). Similar to
 the existence and unique theorem of the stochastic event-triggered
 controlled systems (see, e.g., \cite[Theorem 1]{liuyg2019}),  a unique global solution $(\varrho(t),z(t))$ to the equivalent closed-loop system \dref{closed-loopsystem}
 exists, and then a unique global solution $(x(t),\hat{x}(t))$ to the closed-loop of system \dref{system1.2} under the event-triggered ADRC controller \dref{triggeredadrc}
 also exists.
 Set
\begin{eqnarray}
W_{1}(\varrho)=\varrho Q_{1}\varrho^{\top}, W_{2}(z)=zQ_{2}z^{\top},
\end{eqnarray}
for all $\varrho\in \mathbb{R}^{n}$ and $z\in \mathbb{R}^{n+1}$, where $Q_{1}$ and $Q_{2}$
are the positive definite matrices specified in \dref{Q1E}  and \dref{Q2E}, and  a Lyapunov functional is defined as
\begin{eqnarray}
&&W(t)\triangleq W_{1}(\varrho(t))+W_{2}(z(t))+\int^{t}_{t-\upsilon}\int^{t}_{s}\|\varrho(\sigma)\|^{2}d\sigma ds\cr&&
\int^{t}_{t-\tau}\int^{t}_{s}\|\varrho(\sigma)\|^{2}d\sigma ds
+\int^{t}_{t-\upsilon-\tau}\int^{t}_{s}\|\varrho(\sigma)\|^{2}d\sigma ds\cr&&
+\int^{t}_{t-\upsilon}\int^{t}_{s}\|z(\sigma)\|^{2}d\sigma ds,
\end{eqnarray}
where $\varrho(t)\triangleq\varrho(0),\; t\in [-\upsilon-\tau,0], z(t)\triangleq z(0), t\in [-\upsilon,0].$
Apply It\^{o}'s formula to $W(t)$ with regard to
$t$ along the equivalent closed-loop system \dref{closed-loopsystem} to obtain
\begin{eqnarray*}
&&\hspace{-0.2cm} W(t)
=W(0)-\theta\int^{t}_{0}\|\varrho(s)\|^{2}ds \cr&&\hspace{-0.2cm}+\sum^{n}_{i=1}\int^{t}_{0}\frac{\partial W_{1}(\varrho(s))}{\partial \varrho_{i}}\theta^{n-i}g_{i}(\frac{\varrho_{1}(s)}{\theta^{n-1}},\cdots,\frac{\varrho_{i}(s)}{\theta^{n-i}})ds
\cr&&\hspace{-0.2cm}-\int^{t}_{0}\frac{\partial W_{1}(\varrho(s))}{\partial \varrho_{n}}\sum^{n+1}_{i=1}\frac{\theta^{n+1-i}c_{i}z_{i}(s)}{r^{n+1-i}}ds-\int^{t}_{0}\frac{\partial W_{1}(\varrho(s))}{\partial \varrho_{n}}\cdot\cr&&\sum^{n+1}_{i=1}\theta^{n+1-i}c_{i}\delta_{i}(s)ds
 -r\int^{t}_{0}\|z(s)\|^{2}ds
-\int^{t}_{0}\sum^{n+1}_{i=1}\cr&&\hspace{-0.2cm}\frac{\partial W_{2}(z(s))}{\partial z_{i}}\lambda_{i}r\varsigma(s)ds
+\int^{t}_{0}\sum^{n}_{i=1}\frac{\partial W_{2}(z(s))}{\partial z_{i}}r^{n+1-i}\Xi_{i}(s)ds
\cr&&\hspace{-0.2cm}
+\int^{t}_{0}\frac{\partial W_{2}(z(s))}{\partial z_{n+1}}\Gamma_{1}(s)ds
+\frac{1}{2}\int^{t}_{0}\frac{\partial^{2} W_{2}(z(s))}{\partial z^{2}_{n+1}}[\Gamma^{2}_{2}(s)\cr&&\hspace{-0.2cm}+\Gamma^{2}_{3}(s)]ds+\int^{t}_{0}\frac{\partial W_{2}(z(s))}{\partial z_{n+1}}\Gamma_{2}(s)dB_{1}(s)
+\int^{t}_{0}\frac{\partial W_{2}(z(s))}{\partial z_{n+1}}\cr&&\cdot\Gamma_{3}(s)dB_{2}(s)
-\int^{t}_{0}\int^{s}_{s-v}\|\varrho(\sigma)\|^{2}d\sigma ds+v\int^{t}_{0}\|\varrho(s)\|^{2}ds \cr&&\hspace{-0.2cm}
-\int^{t}_{0}\int^{s}_{s-\tau}\|\varrho(\sigma)\|^{2}d\sigma ds+\tau\int^{t}_{0}\|\varrho(s)\|^{2}ds
\cr&&\hspace{-0.2cm}
-\int^{t}_{0}\int^{s}_{s-\upsilon-\tau}\|\varrho(\sigma)\|^{2}d\sigma ds+(\upsilon+\tau)\int^{t}_{0}\|\varrho(s)\|^{2}ds \cr&&\hspace{-0.2cm}
-\int^{t}_{0}\int^{s}_{s-\upsilon}\|z(\sigma)\|^{2}d\sigma ds+v\int^{t}_{0}\|z(s)\|^{2}ds .
\end{eqnarray*}
Choose sufficiently small $\mu_{1}$, $\mu_{2}$, $\mu_{3}$, $\mu_{4}$ and sufficiently large $r_{1}>0$ such that
\begin{eqnarray}
&& \gamma_{1}\triangleq\theta-2\lambda_{\max}(Q_{1})\sum^{n}_{i=1}L_{i}-\mu_{1}\lambda^{2}_{\max}(Q_{1})-\mu_{2}\lambda^{2}_{\max}(Q_{1})\cr&&-\mu_{4}\lambda^{2}_{\max}(Q_{2})\beta_{2}
-2\upsilon_{1}-2\tau_{1}>0, \cr&&
\gamma_{2}\triangleq1-\mu_{3}\lambda^{2}_{\max}(Q_{2})(\sum^{n+1}_{i=1}|\lambda_{i}|)^{2}>0, \cr&&
\frac{\gamma_{2}r_{1}}{2}-\frac{(\sum^{n+1}_{i=1}\theta^{n+1-i}c_{i})^{2}}{\mu_{1}}-2\lambda_{\max}(Q_{2})\sum^{n}_{i=1}L_{i}\cr&&-\mu_{4}\lambda^{2}_{\max}(Q_{2})\beta_{3}-\frac{1}{\mu_{4}}-\upsilon_{1}>0,
\end{eqnarray}
where $\tau_{1}\triangleq \epsilon_{1}r^{-(2n+\frac{3}{2})}_{1}$ and $\upsilon_{1}\triangleq \epsilon_{2}r^{-(\frac{2n}{3}+\frac{5}{3})}_{1}$.
These together with \dref{55cite}, \dref{56cite}, and Young's inequality, further yield that for all $r\geq r_{1}$,
\begin{eqnarray}\label{fe1331}
&&\hspace{-0.2cm} \frac{d\mathbb{E}W(t)}{dt}
\cr&&\hspace{-0.2cm} \leq -\theta\mathbb{E}\|\varrho(t)\|^{2}+
2\lambda_{\max}(Q_{1})\sum^{n}_{i=1}L_{i}\mathbb{E}\|\varrho(t)\|^{2}
\cr&&\hspace{-0.2cm} +\mu_{1}\lambda^{2}_{\max}(Q_{1})\mathbb{E}\|\varrho(t)\|^{2}
+\frac{(\sum^{n+1}_{i=1}\theta^{n+1-i}c_{i})^{2}}{\mu_{1}}\mathbb{E}\|z(t)\|^{2}\cr&&\hspace{-0.2cm}
+\mu_{2}\lambda^{2}_{\max}(Q_{1})\mathbb{E}\|\varrho(t)\|^{2}+\frac{1}{\mu_{2}}\mathbb{E}[\sum^{n+1}_{i=1}|\theta^{n+1-i}c_{i}\delta_{i}(t)|\mathbb{I}_{\Omega^{*}_{\upsilon}}]^{2}\cr&&\hspace{-0.2cm}
+\mathbb{E}[\sum^{n+1}_{i=1}|\frac{\theta^{n+1-i}c_{i}}{\mu_{2}}\delta_{i}(t)|\mathbb{I}_{\Omega\setminus\Omega^{*}_{\upsilon}}]^{2}
 -r\mathbb{E}\|z(t)\|^{2}+\mu_{3}\lambda^{2}_{\max}(Q_{2})\cr&&\hspace{-0.2cm}\cdot(\sum_{i=1}^{n+1}|\lambda_{i}|)^{2}r\mathbb{E}\|z(t)\|^{2}
 +\frac{r}{\mu_{3}}\mathbb{E}[\varsigma^{2}(t)\mathbb{I}_{\Omega_{\tau}}]
 +\frac{r}{\mu_{3}}\mathbb{E}[\varsigma^{2}(t)\mathbb{I}_{\Omega\setminus\Omega_{\tau}}]
 \cr&&\hspace{-0.2cm}+2\lambda_{\max}(Q_{2})\sum^{n}_{i=1}L_{i}\mathbb{E}\|z(t)\|^{2}
 +\mu_{4}\lambda^{2}_{\max}(Q_{2})[\beta_{1}+\beta_{2}\mathbb{E}\|\varrho(t)\|^{2}\cr&&\hspace{-0.2cm}
 +\beta_{3}\mathbb{E}\|z(t)\|^{2}+\beta_{4}\mathbb{E}\|\delta(t)\|^{2}]+\frac{1}{\mu_{4}}\mathbb{E}\|z(t)\|^{2}+\lambda_{\max}(Q_{2})\beta_{5}\cr&&\hspace{-0.2cm}-\mathbb{E}\int^{t}_{t-\upsilon}\|\varrho(s)\|^{2}ds
 -\mathbb{E}\int^{t}_{t-\tau}\|\varrho(s)\|^{2}ds
 +2(\upsilon +\tau)\mathbb{E}\|\varrho(t)\|^{2}\cr&&\hspace{-0.2cm}
-\mathbb{E}\int^{t}_{t-\upsilon-\tau}\|\varrho(s)\|^{2}ds
-\mathbb{E}\int^{t}_{t-\upsilon}\|z(s)\|^{2}ds+\upsilon\mathbb{E}\|z(t)\|^{2}\cr&&\hspace{-0.2cm}
\leq -\gamma_{1}\mathbb{E}\|\varrho(t)\|^{2}-\frac{\gamma_{2}r}{2}\mathbb{E}\|z(t)\|^{2}-\mathbb{E}\int^{t}_{t-\upsilon}\|\varrho(s)\|^{2}ds\cr&&\hspace{-0.2cm}-\mathbb{E}\int^{t}_{t-\tau}\|\varrho(s)\|^{2}ds
-\mathbb{E}\int^{t}_{t-\upsilon-\tau}\|\varrho(s)\|^{2}ds\cr&&\hspace{-0.2cm}-\mathbb{E}\int^{t}_{t-\upsilon}\|z(s)\|^{2}ds
+\frac{1}{\mu_{2}}\mathbb{E}[\sum^{n+1}_{i=1}|\theta^{n+1-i}c_{i}\delta_{i}(t)|\mathbb{I}_{\Omega^{*}_{\upsilon}}]^{2}
\cr&&\hspace{-0.2cm}+\frac{1}{\mu_{2}}\mathbb{E}[\sum^{n+1}_{i=1}|\theta^{n+1-i}c_{i}\delta_{i}(t)|\cdot\mathbb{I}_{\Omega\setminus\Omega^{*}_{\upsilon}}]^{2}+
 \frac{r}{\mu_{3}}\mathbb{E}[\varsigma^{2}(t)\mathbb{I}_{\Omega_{\tau}}]
\cr&&\hspace{-0.2cm}+\frac{r}{\mu_{3}}\mathbb{E}[\varsigma^{2}(t)\mathbb{I}_{\Omega\setminus\Omega_{\tau}}]+\mu_{4}\lambda^{2}_{\max}(Q_{2})\beta_{4}\cdot\mathbb{E}[\|\delta(t)\|^{2}\mathbb{I}_{\Omega^{*}_{\upsilon}}]
\cr&&\hspace{-0.2cm}+\mu_{4}\lambda^{2}_{\max}(Q_{2}) \beta_{4}\mathbb{E}[\|\delta(t)\|^{2}\mathbb{I}_{\Omega\setminus\Omega^{*}_{\upsilon}}]
+\mu_{4}\lambda^{2}_{\max}(Q_{2})\beta_{1}\cr&&\hspace{-0.2cm}+\lambda_{\max}(Q_{2})\beta_{5}.
\end{eqnarray}
By the ETMs \dref{triggermechemi} and \dref{triggermechem2}, it follows that for all $t\geq 0$,
\begin{eqnarray}\label{dfe332}
&&  \mathbb{E}[\varsigma^{2}(t)\mathbb{I}_{\Omega\setminus\Omega_{\tau}}]=\sum^{[\frac{t}{\tau}]+1}_{k=1}\mathbb{E}[\varsigma^{2}(t)\mathbb{I}_{\Omega_{k}\setminus\Omega_{k,\tau}}]\leq \frac{\kappa^{2}_{1}}{r},\;
\cr && \mathbb{E}[\sum^{n+1}_{i=1}|\delta_{i}(t)|\mathbb{I}_{\Omega\setminus\Omega^{*}_{v}}]^{2}
\leq \frac{\kappa^{2}_{2}}{r}.
\end{eqnarray}
For $\omega\in \Omega_{k,\tau}$, we have $t_{k}\geq t-\tau$. A direct computation shows that
\begin{eqnarray}\label{dfdfe2}
&& \mathbb{E}[\varsigma^{2}(t)\mathbb{I}_{\Omega_{\tau}}]=\sum^{[\frac{t}{\tau}]+1}_{k=1}\mathbb{E}[\varsigma^{2}(t)\mathbb{I}_{\Omega_{k,\tau}}]\cr &&
=r^{2n}\sum^{[\frac{t}{\tau}]+1}_{k=1}\mathbb{E}\left[\left(\int^{t}_{t_{k}}x_{2}(s)ds\right)^{2}\mathbb{I}_{\Omega_{k,\tau}}\right] \cr&&
\leq r^{2n}\tau \int^{t}_{t-\tau}\mathbb{E}\|\varrho(s)\|^{2}ds.
\end{eqnarray}
Set
\begin{eqnarray}\label{515equaitonsf}
\gamma^{*}=\frac{\theta^{2n}}{\mu_{2}}\disp\max_{1\leq i\leq n+1}c^{2}_{i}+\mu_{4}\lambda^{2}_{\max}(Q_{2})\beta_{4}.
\end{eqnarray}

By \dref{414df}, $\tau=\epsilon_{1} r^{-(2n+\frac{3}{2})}$, and $\upsilon=\epsilon_{2}r^{-(\frac{2n}{3}+\frac{5}{3})}$, it can be obtained that
$\Lambda_{1}(\upsilon)$, $\Lambda_{2}(\upsilon)$, $\Lambda_{3}(\upsilon,\tau,r)$, and $\Lambda_{5}(\upsilon,r)$ are strictly decreasing
with respect to $r$ and approach zero as $r\rightarrow \infty$.
Therefore, there exists  $r_{2}>0$ such that
\begin{eqnarray}\label{432df}
&&\gamma_{3}\triangleq\gamma_{1}-\gamma^{*}\Lambda_{1}(\upsilon_{2})>0,
\gamma_{4}\triangleq1-\gamma^{*}\Lambda_{2}(\upsilon_{2})>0,\cr&&
\gamma_{5}\triangleq1-\frac{\epsilon_{1}}{\mu_{3}r^{\frac{1}{2}}_{2}}>0,
\gamma_{6}\triangleq1-\gamma^{*}\Lambda_{3}(\upsilon_{2},\tau_{2},r_{2})>0, \cr&&
\gamma_{7}\triangleq1-\gamma^{*}\Lambda_{5}(\upsilon_{2},r_{2})>0,
\end{eqnarray}
where $\tau_{2}\triangleq \epsilon_{1}r^{-(2n+\frac{3}{2})}_{2}$, $\upsilon_{2}\triangleq \epsilon_{2}r^{-(\frac{2n}{3}+\frac{5}{3})}_{2}$,
and $r_{2}$ is dependent on $\gamma_{1}$, $\gamma^{*}$, $n$, $\theta$, $\alpha_{2}$, $\epsilon_{1}$, $\epsilon_{2}$, $\mu_{3}$, $L_{i}\;(i=1,\cdots,n)$, $c_{i}\;(i=1,\cdots,n+1)$, and $\lambda_{i}\;(i=1,\cdots,n+1)$ .
Choose
\begin{eqnarray}\label{516equation}
r\geq r_{*}\triangleq \max\{1,\frac{4\gamma^{*}}{\gamma_{2}}\Lambda_{4}(\upsilon_{2}),\frac{\epsilon^{2}_{1}}{\mu^{2}_{3}},r_{1},r_{2}\}.
\end{eqnarray}
It follows from Lemma \ref{errorlemma}, \dref{fe1331}, \dref{dfe332}, \dref{dfdfe2}, \dref{515equaitonsf},  \dref{432df}, and \dref{516equation} that
\begin{eqnarray}\label{335dd}
&&\hspace{-0.2cm} \frac{d\mathbb{E}W(t)}{dt}
\cr&&\hspace{-0.2cm}\leq -\gamma_{1}\mathbb{E}\|\varrho(t)\|^{2}-\frac{\gamma_{2}r}{2}\mathbb{E}\|z(t)\|^{2}-\mathbb{E}\int^{t}_{t-\upsilon}\|\varrho(s)\|^{2}ds\cr&&\hspace{-0.2cm}-\mathbb{E}\int^{t}_{t-\tau}\|\varrho(s)\|^{2}ds
-\mathbb{E}\int^{t}_{t-\upsilon-\tau}\|\varrho(s)\|^{2}ds-\mathbb{E}\int^{t}_{t-\upsilon}\|z(s)\|^{2}\cr&&\hspace{-0.2cm}ds+\gamma^{*}\bigg\{\Lambda_{1}(\upsilon_{2})\mathbb{E}[\|\varrho(t)\|^{2}\mathbb{I}_{\Omega^{*}_{\upsilon}}]
+\Lambda_{2}(\upsilon_{2})\int^{t}_{t-\upsilon}\mathbb{E}[\|\varrho(s)\|^{2}\cr&&\hspace{-0.2cm}\mathbb{I}_{\Omega^{*}_{\upsilon}}]ds
+\Lambda_{3}(\upsilon_{2},\tau_{2},r_{2})\int^{t}_{t-\upsilon-\tau}\mathbb{E}[\|\varrho(s)\|^{2}\mathbb{I}_{\Omega^{*}_{\upsilon}}]ds+\Lambda_{4}(\upsilon_{2})\cdot\cr&&\hspace{-0.2cm}\mathbb{E}[\|z(t)\|^{2}\mathbb{I}_{\Omega^{*}_{\upsilon}}]
+\Lambda_{5}(\upsilon_{2},r_{2})\int^{t}_{t-\upsilon}\mathbb{E}[\|z(s)\|^{2}\mathbb{I}_{\Omega^{*}_{\upsilon}}]ds+
\Lambda_{6}(\upsilon_{2},\cr&&\hspace{-0.2cm}r_{2})\bigg\}+\max_{1\leq i\leq n+1}c^{2}_{i}\frac{\theta^{2n}\kappa^{2}_{2}}{\mu_{2}r}
+\frac{r^{2n+1}\tau}{\mu_{3}}\int^{t}_{t-\tau}\mathbb{E}\|\varrho(s)\|^{2}ds
 +\frac{\kappa^{2}_{1}}{\mu_{3}} \cr&&\hspace{-0.2cm}+\frac{\mu_{4}\lambda^{2}_{\max}(Q_{2}) \beta_{4}\kappa^{2}_{2}}{r}
 +\mu_{4}\lambda^{2}_{\max}(Q_{2})\beta_{1}+\lambda_{\max}(Q_{2})\beta_{5}\cr&&\hspace{-0.2cm}
 \leq -\gamma_{3}\mathbb{E}\|\varrho(t)\|^{2}-\frac{\gamma_{2}r}{4}\mathbb{E}\|z(t)\|^{2}-\gamma_{4}\mathbb{E}\int^{t}_{t-\upsilon}\|\varrho(s)\|^{2}ds\cr&&\hspace{-0.2cm}-\gamma_{5}\mathbb{E}\int^{t}_{t-\tau}\|\varrho(s)\|^{2}ds
 -\gamma_{6}\mathbb{E}\int^{t}_{t-\upsilon-\tau}[\|\varrho(s)\|^{2}ds\cr&&\hspace{-0.2cm}-\gamma_{7}\mathbb{E}\int^{t}_{t-\upsilon}\|z(s)\|^{2}ds+M_{1}\cr&&\hspace{-0.2cm}\leq
-\frac{\gamma_{3}}{\lambda_{\max}(Q_{1})}\mathbb{E}W_{1}(\varrho(t))-\frac{\gamma_{2}r}{4\lambda_{\max}(Q_{2})}\mathbb{E}W_{2}(z(t))
\cr&&\hspace{-0.2cm}-\frac{\gamma_{4}}{\upsilon}\mathbb{E}\int^{t}_{t-\upsilon}\int^{t}_{s}\|\varrho(\sigma)\|^{2}d\sigma ds-\frac{\gamma_{5}}{\tau}\mathbb{E}\int^{t}_{t-\tau}\int^{t}_{s}\|\varrho(\sigma)\|^{2}d\sigma ds
\cr&&\hspace{-0.2cm}-\frac{\gamma_{6}}{\upsilon+\tau} \mathbb{E}\int^{t}_{t-\upsilon-\tau}\int^{t}_{s}\|\varrho(\sigma)\|^{2}d\sigma ds
\cr&&\hspace{-0.2cm} -\frac{\gamma_{7}}{\upsilon}\mathbb{E}\int^{t}_{t-\upsilon}\int^{t}_{s}\|z(\sigma)\|^{2}d\sigma ds+M_{1}
\cr&&\hspace{-0.2cm} \leq -\gamma\mathbb{E}W(t)+M_{1},
\end{eqnarray}
where we set
\begin{eqnarray*}
&& \gamma=\min\{\frac{\gamma_{3}}{\lambda_{\max}(Q_{1})},\frac{\gamma_{2}r_{*}}{4\lambda_{\max}(Q_{2})},\gamma_{4},\gamma_{5},\frac{\gamma_{6}}{2},\gamma_{7}\}, \cr
&&M_{1}=\gamma^{*}\Lambda_{6}(\upsilon_{2},r_{2})+\max_{1\leq i\leq n+1}c^{2}_{i}\frac{\theta^{2n}\kappa^{2}_{2}}{\mu_{2}r_{*}}+\frac{\kappa^{2}_{1}}{\mu_{3}}\cr&&+\frac{\mu_{4}\lambda^{2}_{\max}(Q_{2}) \beta_{4}\kappa^{2}_{2}}{r_{*}}
 +\mu_{4}\lambda^{2}_{\max}(Q_{2})\beta_{1}+\lambda_{\max}(Q_{2})\beta_{5}.
\end{eqnarray*}
%
For any fixed $\varepsilon>0$, a direct computation shows that
\begin{eqnarray}\label{337df}
&&\sup_{r\geq r_{*}}e^{-\gamma r^{\varepsilon}}\mathbb{E}W(0)\cr&&\leq \sup_{r\geq r_{*}}e^{-\gamma r^{\varepsilon}}\bigg\{[\lambda_{\max}(Q_{1})+\upsilon^{2}+\tau^{2}+(\upsilon+\tau)^{2}]\|\varrho(0)\|^{2}\cr&&+(\lambda_{\max}(Q_{2})+\upsilon^{2})\sum^{n+1}_{i=1}r^{2(n +1-i)}\mathbb{E}|x_{i}(0)-\hat{x}_{i}(0)|^{2}\bigg\}
\cr&&\triangleq M_{2},
\end{eqnarray}
where $M_{2}<\infty$ because $\disp\lim_{r\rightarrow \infty}e^{-\gamma r^{\varepsilon}}r^{2n}=0$. Thus, for all $t\geq r^{\varepsilon}$ with $r\geq r_{*}$, it follows from \dref{335dd} and \dref{337df} that
\begin{eqnarray}\label{boundedness436}
&& \mathbb{E}W(t)\leq e^{-\gamma t}\mathbb{E}W(0)+\int^{t}_{0} e^{-\gamma (t-s)}M_{1}ds
\cr&&\leq  M_{2}+\frac{M_{1}}{\gamma}\triangleq M_{3}.
\end{eqnarray}
By Lemma \ref{errorlemma}, ETM \dref{triggermechem2}, and \dref{boundedness436},
for all $t\geq r^{\varepsilon}+2\geq r^{\varepsilon}+\upsilon_{2}+\tau_{2}$, we have
\begin{eqnarray}\label{437bounded}
&&\hspace{-0.1cm}\mathbb{E}\|\delta(t)\|^{2}=\mathbb{E}[\|\delta(t)\|^{2}\mathbb{I}_{\Omega^{*}_{\upsilon}}]+\mathbb{E}[\|\delta(t)\|^{2}\mathbb{I}_{\Omega\setminus\Omega^{*}_{\upsilon}}]\cr&&\hspace{-0.1cm}
\leq \bigg[\frac{\Lambda_{1}(\upsilon_{2})}{\lambda_{\min}(Q_{1})}+\frac{\Lambda_{2}(\upsilon_{2})\upsilon }{\lambda_{\min}(Q_{1})}
+\frac{\Lambda_{3}(\upsilon_{2},\tau_{2},r_{2})(\upsilon_{2}+\tau_{2}) }{\lambda_{\min}(Q_{1})}\cr&&\hspace{-0.1cm}+\frac{\Lambda_{4}(\upsilon_{2}) }{\lambda_{\min}(Q_{2})}+\frac{\Lambda_{5}(\upsilon_{2},r_{2})\upsilon_{2} }{\lambda_{\min}(Q_{2})}\bigg]M_{3}
+\Lambda_{6}(\upsilon_{2},r_{2})+\frac{\kappa^{2}_{2}}{r^{2}_{*}}\cr&&\hspace{-0.1cm}\triangleq M_{4}.
\end{eqnarray}
Similar to \dref{fe1331}, it follows from \dref{dfe332}, \dref{dfdfe2}, \dref{boundedness436}, and \dref{437bounded} that for all $t\geq r^{\varepsilon}+2$,
\begin{eqnarray}
&&\hspace{-0.2cm}\frac{d\mathbb{E}W_{2}(z(t))}{dt}
\cr&&\hspace{-0.2cm}\leq -(1-\mu_{3}\lambda^{2}_{\max}(Q_{2})(\sum^{n+1}_{i=1}|\lambda_{i}|)^{2})r\mathbb{E}\|z(t)\|^{2}
 \cr&&\hspace{-0.2cm}+\frac{r}{\mu_{3}}\mathbb{E}[\varsigma^{2}(t)\mathbb{I}_{\Omega_{\tau}}]
 +\frac{r}{\mu_{3}}\mathbb{E}[\varsigma^{2}(t)\mathbb{I}_{\Omega\setminus\Omega_{\tau}}]+2\lambda_{\max}(Q_{2})\cdot\cr&&\hspace{-0.2cm}\sum^{n}_{i=1}L_{i}\mathbb{E}\|z(t)\|^{2} +\mu_{4}\lambda^{2}_{\max}(Q_{2})[\beta_{1}
 +\beta_{2}\mathbb{E}\|\varrho(t)\|^{2}
+\beta_{3}\cdot\cr&&\hspace{-0.2cm}\mathbb{E}\|z(t)\|^{2}+\beta_{4}\mathbb{E}\|\delta(t)\|^{2}]+\frac{1}{\mu_{4}}\mathbb{E}\|z(t)\|^{2}+\lambda_{\max}(Q_{2})\beta_{5}
 \cr&&\hspace{-0.2cm}\leq -\gamma_{2}r\mathbb{E}\|z(t)\|^{2}+\frac{\epsilon_{1}}{r^{\frac{1}{2}}_{*}\mu_{3}}\int^{t}_{t-\tau}\mathbb{E}\|\varrho(s)\|^{2}ds+\frac{\kappa^{2}_{1}}{\mu_{3}}\cr&&\hspace{-0.2cm}+2\lambda_{\max}(Q_{2})\sum^{n}_{i=1}L_{i}\mathbb{E}\|z(t)\|^{2}
 +\mu_{4}\lambda^{2}_{\max}(Q_{2})[\beta_{1}+\beta_{2}\mathbb{E}\|\varrho(t)\|^{2}
\cr&&\hspace{-0.2cm} +\beta_{3}\mathbb{E}\|z(t)\|^{2}+\beta_{4}\mathbb{E}\|\delta(t)\|^{2}]+\frac{1}{\mu_{4}}\mathbb{E}\|z(t)\|^{2}+\lambda_{\max}(Q_{2})\beta_{5}\cr&&\hspace{-0.2cm}
 \leq  -\frac{\gamma_{2}r}{\lambda_{\max}(Q_{2})}\mathbb{E}W_{2}(z(t))+M_{5},
 \end{eqnarray}
 where we set
 \begin{eqnarray}
 &&\hspace{-0.2cm}M_{5}\triangleq\frac{\epsilon_{1}M_{3}}{r^{\frac{1}{2}}_{*}\mu_{3}\lambda_{\min}(Q_{1})}+\frac{\kappa^{2}_{1}}{\mu_{3}}+\frac{2\lambda_{\max}(Q_{2})\sum^{n}_{i=1}L_{i}M_{3}}{\lambda_{\min}(Q_{2})}
 \cr&&\hspace{-0.2cm}+u_{4}\lambda^{2}_{\max}(Q_{2})[\beta_{1}+\frac{\beta_{2}M_{3}}{\lambda_{\min}(Q_{1})}
+\frac{\beta_{3}M_{3}}{\lambda_{\min}(Q_{2})}+\beta_{4}M_{4}] \cr&&\hspace{-0.2cm} +\frac{M_{3}}{\mu_{4}\lambda_{\min}(Q_{2})}+\lambda_{\max}(Q_{2})\beta_{5}.
 \end{eqnarray}
This together with \dref{boundedness436}, yields that for all  $t\geq r^{\varepsilon}+3$,
 \begin{eqnarray}\label{343equation}
&&\hspace{-0.2cm} \mathbb{E}W_{2}(z(t))\cr&&\hspace{-0.2cm}\leq e^{-\frac{\gamma_{2}r(t-r^{\varepsilon}-2)}{\lambda_{\max}(Q_{2})}}\mathbb{E}W_{2}(z(r^{\varepsilon}+2))
+\int^{t}_{r^{\varepsilon}+2}e^{-\frac{\gamma_{2}r(t-s)}{\lambda_{\max}(Q_{2})}}M_{5}ds\cr&&\hspace{-0.2cm}
\leq \frac{1}{r}[M_{6}+\frac{\lambda_{\max}(Q_{2})M_{5}}{\gamma_{2}}]\triangleq\frac{M_{7}}{r},
 \end{eqnarray}
where $M_{6}\triangleq\sup_{r\geq r_{*}}re^{-\frac{\gamma_{2}r}{\lambda_{\max}(Q_{2})}}M_{3}$. Therefore, for all $i=1,\cdots,n+1$ and
$t\geq r^{\varepsilon}+3$,
\begin{eqnarray}
&& \hspace{-0.2cm}\mathbb{E}|x_{i}(t)-\hat{x}_{i}(t)|^{2}=\frac{1}{r^{2(n+1-i)}}\mathbb{E}|z_{i}(t)|^{2}\cr&&\hspace{-0.2cm}\leq \frac{\mathbb{E}W_{2}(z(t))}{r^{2(n+1-i)}\lambda_{\min}(Q_{2})}
\leq \frac{M_{7}}{\lambda_{\min}(Q_{2})r^{2n+3-2i}}.
 \end{eqnarray}
 Similar to \dref{fe1331},  for all $t\geq r^{\varepsilon}+5$, it follows from Lemma \ref{errorlemma}, \dref{dfe332}, \dref{boundedness436}, and \dref{343equation} that
 \begin{eqnarray}\label{344fd}
&&\hspace{-0.2cm}\frac{d\mathbb{E}W_{1}(\varrho(t))}{dt}\cr&&\hspace{-0.2cm}
\leq -\gamma_{1}\mathbb{E}\|\varrho(t)\|^{2}
+\frac{(\sum^{n+1}_{i=1}\theta^{n+1-i}c_{i})^{2}}{\mu_{1}}\mathbb{E}\|z(t)\|^{2}
+\frac{1}{\mu_{2}}\mathbb{E}[\sum^{n+1}_{i=1}\cr&&\hspace{-0.2cm}|
\theta^{n+1-i}c_{i}\delta_{i}(t)|\mathbb{I}_{\Omega^{*}_{v}}]^{2}
+\frac{1}{\mu_{2}}\mathbb{E}[\sum^{n+1}_{i=1}
|\theta^{n+1-i}c_{i}\delta_{i}(t)|\mathbb{I}_{\Omega\setminus\Omega^{*}_{v}}]^{2}
\cr&&\hspace{-0.2cm}\leq -\frac{\gamma_{1}}{\lambda_{\max}(Q_{1})}\mathbb{E}W_{1}(\varrho(t))+\frac{M_{8}}{r},
  \end{eqnarray}
where
\begin{eqnarray}
&&\hspace{-0.2cm}M_{8}\triangleq\frac{(\sum^{n+1}_{i=1}\theta^{n+1-i}c_{i})^{2}M_{7}}{\mu_{1}\lambda_{\min}(Q_{2})}+\frac{\theta^{2n}\disp\max_{1\leq i\leq n+1}c^{2}_{i}}{\mu_{2}}\sup_{r\geq r_{*}}\cr&&\hspace{-0.2cm}\bigg\{
\frac{\Lambda_{1}(\upsilon)rM_{3}}{\lambda_{\min}(Q_{1})}+
\frac{\Lambda_{2}(\upsilon)vrM_{3}}{\lambda_{\min}(Q_{1})}
+\frac{\Lambda_{3}(\upsilon,\tau,r)(v+\tau)rM_{3}}{\lambda_{\min}(Q_{1})}\cr&&\hspace{-0.2cm}+\frac{\Lambda_{4}(\upsilon)M_{7}}{\lambda_{\min}(Q_{2})}+\frac{\Lambda_{5}(\upsilon,r)\upsilon M_{7}}{\lambda_{\min}(Q_{2})}+\Lambda_{6}(\upsilon,r)r+\kappa^{2}_{2}\bigg\}<\infty
\end{eqnarray}
by \dref{414df} and $\upsilon=\epsilon_{2}r^{-(\frac{2n}{3}+\frac{5}{3})}.$
Thus, it follows from \dref{boundedness436} and \dref{344fd} that for all
\begin{eqnarray}t\geq t_{r}\triangleq 2r^{\varepsilon}+5,\end{eqnarray} with $r\geq r_{*}$ and $\varepsilon$ be
any positive constant, we have
 \begin{eqnarray}
&&\hspace{-0.2cm} \mathbb{E}W_{1}(\varrho(t))\cr&&\hspace{-0.2cm} \leq e^{-\frac{\gamma_{1}}{\lambda_{\max}(Q_{1})}(t-r^{\varepsilon}-5)}W_{1}(\varrho(r^{\varepsilon}+5))
+\int^{t}_{r^{\varepsilon}+5}e^{-\frac{\gamma_{1}(t-s)}{\lambda_{\max}(Q_{1})}}\cr&&\hspace{-0.2cm}\frac{M_{8}}{r}ds
\leq e^{-\frac{\gamma_{1}}{\lambda_{\max}(Q_{1})}r^{\varepsilon}}M_{3}+\frac{\lambda_{\max}(Q_{1})M_{8}}{\gamma_{1}r}\leq \frac{M_{9}}{r},
  \end{eqnarray}
where
 \begin{eqnarray}
M_{9}\triangleq \sup_{r\geq r^{*}}re^{-\frac{\gamma_{1}}{\lambda_{\max}(Q_{1})}r^{\varepsilon}}M_{3}+\frac{\lambda_{\max}(Q_{1})M_{8}}{\gamma_{1}}.
  \end{eqnarray}
Thus, for all $t\geq t_{r}$ with $r\geq r_{*}$, it holds that
 \begin{eqnarray}
\sum^{n}_{i=1}\mathbb{E}|x_{i}(t)|^{2}\leq \sum^{n}_{i=1} \mathbb{E}|\varrho_{i}(t)|^{2} \leq \frac{\mathbb{E}W_{1}(\varrho(t))}{\lambda_{\min}(Q_{1})}\leq \frac{M}{r},
  \end{eqnarray}
where we set
 \begin{eqnarray}
M=\frac{M_{9}}{\lambda_{\min}(Q_{1})},
  \end{eqnarray}
which is independent of the tuning gain $r$.
This finishes the proof of (i) and (iii) of Theorem \ref{theorem34.1}.
It follows from (i) of Theorem \ref{theorem34.1}
 and Chebyshev's inequality (\cite[p.5]{mao}) that
\begin{eqnarray}
P\{|x_{i}(t)-\hat{x}_{i}(t)|\geq k^{2}\sqrt{M}r^{-\frac{2n+3-2i}{2}}\} \leq \frac{1}{k^{2}}
\end{eqnarray}
for all $t\geq t_{r}$, $k\in \mathbb{Z}^{+}$, and $i=1,\cdots,n+1$.  From the Borel-Cantelli's lemma (\cite[p.7]{mao}), for almost all $\omega\in \Omega$
 and $t\geq t_{r}$, there
is a random variable $k_{0}(\omega)$ such that whenever $k\geq k_{0}(\omega)$, we have
\begin{eqnarray}
|x_{i}(t)-\hat{x}_{i}(t)|\leq  k^{2}\sqrt{M}r^{-\frac{2n+3-2i}{2}}.
\end{eqnarray}
Thus,  (ii) of Theorem \ref{theorem34.1} holds and  (iv) follows similarly with $M_{\omega}\triangleq k^{2}_{0}(\omega)\sqrt{M}>0$ be a random variable independent of $r$.
This completes the proof of Theorem \ref{theorem34.1}.
\hfill $\Box$

\section{Numerical simulations}\label{sectionIV}
Some numerical simulations are implemented to verify the functionality of the
proposed event-triggered ADRC scheme in this section. The following uncertain random nonlinear systems are taken as an numerical example:
\begin{equation}\label{numericalexample}
\hspace{-0.2cm}\left\{\begin{array}{l} \dot{x}_{1}(t)=x_{2}(t)+\sin(x_{1}(t)), \cr
\dot{x}_{2}(t)=f(t,x(t),w_{1}(t),w_{2}(t))+\sin(x_{1}(t)+x_{2}(t))\cr\hspace{1.2cm}+u(t),\cr
y(t)=x_{1}(t),
\end{array}\right.
\end{equation}
which is the second-order case of system \dref{system1.2} with $g_{1}(x_{1})=\sin(x_{1})$, $g_{2}(x)=\sin(x_{1}+x_{2})$.
$x_{3}(t)\triangleq f(t,x(t),w_{1}(t),w_{2}(t))$ is the random total disturbance.
Choose $r=50$, $\lambda_{1}=6$, $\lambda_{2}=12$, $\lambda_{3}=8$, and $\epsilon_{1}=\kappa_{1}=1$. Thus, the eigenvalues of $H$ in \dref{matricf2} are equal to $-2$ and then $H$ is Hurwitz.
 The event-triggered ESO is then designed as
 \begin{equation}
\left\{\begin{array}{l}
\dot{\hat{x}}_{1}(t)=\hat{x}_{2}(t)+6{\color{blue}\cdot}50\left(y(t_{k})-\hat{x}_{1}(t)\right)+\sin(\hat{x}_{1}(t)),
\crr
\dot{\hat{x}}_{2}(t)=\hat{x}_{3}(t)+12{\color{blue}\cdot}50^{2}\left(y(t_{k})-\hat{x}_{1}(t)\right)\crr\hspace{1.2cm}+\sin(\hat{x}_{1}(t)+\hat{x}_{2}(t))+u(t),\crr
 \dot{\hat{x}}_{3}(t)=8{\color{blue}\cdot}50^{3}\left(y(t_{k})-\hat{x}_{1}(t)\right),
\end{array}\right.
\end{equation}
where execution times $t_{k}\;(k\in \mathbb{Z}^{+})$ are determined by the following ETM
\begin{equation}\label{ETMESOSimul}
t_{k+1}=\inf\{t\geq t_{k}+50^{-5.5}:\;|y(t)-y(t_{k})|\geq 50^{-2.5} \}.
\end{equation}
Design $c_{1}=-1$, $c_{2}=-2$, and then $Q_{1}=
\begin{pmatrix}
\frac{3}{2}     & \frac{1}{2}                \cr \frac{1}{2}    & \frac{1}{2}
\end{pmatrix}$
with $\lambda_{\max}(Q_{1})=1.7071$. Take $\theta=7$ which satisfies $\theta>2\lambda_{\max}(Q_{1})\sum^{n}_{i=1}L_{i}=6.8284$
and $\upsilon\theta^{2}\max_{1\leq i\leq 2}\{|c_{i}|,1\}=50^{-3}\times7^{2}\times2 \approx 0.000784<1$, and let $\epsilon_{2}=\kappa_{2}=1$.
Therefore, the  event-triggered ADRC controller is designed as
\begin{eqnarray}\label{eventADRCcontroller}
u(t)=-49\hat{x}_{1}(t^{*}_{l})-14\hat{x}_{2}(t^{*}_{l})-\hat{x}_{3}(t^{*}_{l}), \; t\in [t^{*}_{l},t^{*}_{l+1}),
\end{eqnarray}
where
\begin{equation}\label{ETMADRCSimul}
t^{*}_{l+1}=\inf\{t\geq t^{*}_{l}+50^{-3}:\;\sum^{3}_{i=1}|\hat{x}_{i}(t)-\hat{x}_{i}(t^{*}_{l})|\geq \frac{1}{50^{\frac{1}{2}}}\}.
\end{equation}
In the following numerical simulations, we take
\begin{eqnarray}\label{numstoachsa}
&&\hspace{-0.2cm}f(t,x,w_{1},w_{2})=x_{1}+2x_{2}+\sin(t)+\cos(x_{1}+x_{2})\cr &&\hspace{-0.2cm}+w^{3}_{1}+w_{2}, \;\; w_{1}(t)=2\sin(t+B_{1}(t)),
\end{eqnarray}
and $\varrho_{1}=\varrho_{2}=1.5$ for definition of $w_{2}(t)$ in \dref{21equ}, and the initial values are specified as
$x_{1}(0)=0.5,x_{2}(0)=-0.5,w_{2}(0)=0,\hat{x}_{1}(0)=\hat{x}_{2}(0)=\hat{x}_{3}(0)=0$.
It can be easily checked that all assumptions of Theorem \ref{theorem34.1} are satisfied.

\begin{figure}[ht]\centering
\subfigure[]
 {\includegraphics[width=4.1cm,height=4cm]{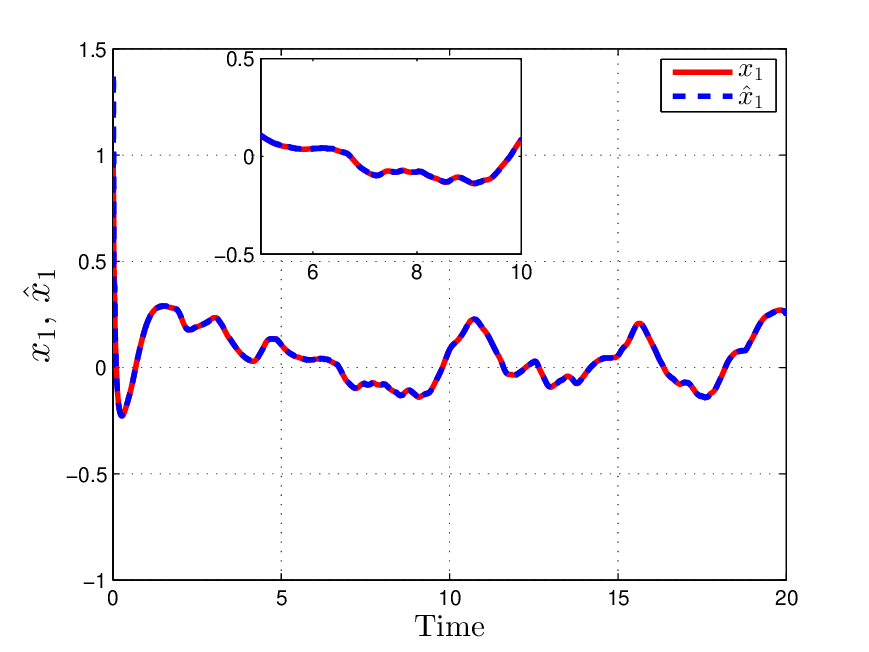}\label{Fig-Lx1}}
\subfigure[]
 {\includegraphics[width=4.1cm,height=4cm]{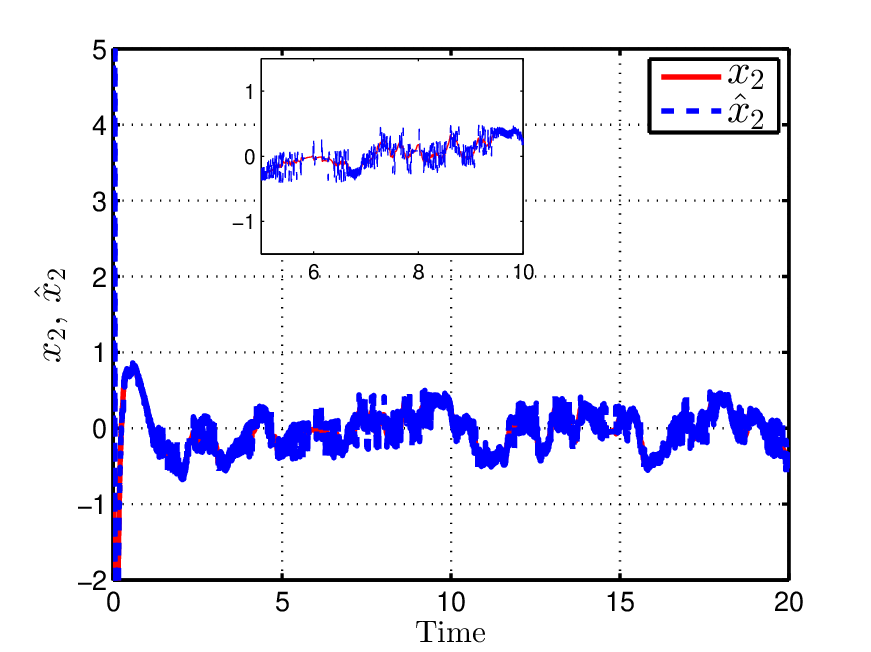}\label{Fig-Lx2}}
 \subfigure[]
 {\includegraphics[width=4.1cm,height=4cm]{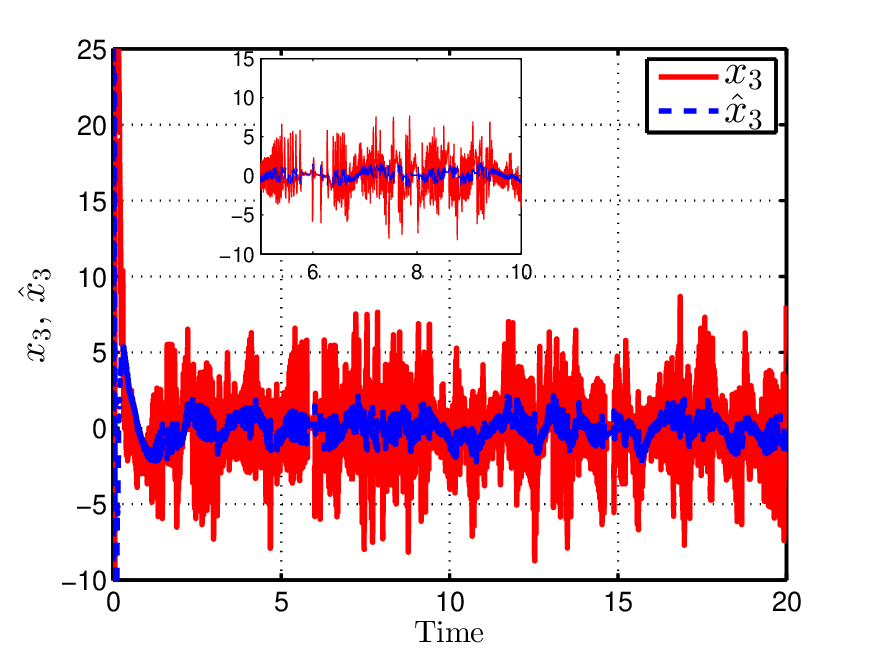}\label{Fig-Lx3}}
\caption{The state and random total disturbance $(x_{1}(t),x_{2}(t),x_{3}(t))$
and their estimates $(\hat{x}_1(t),\hat{x}_{2}(t),\hat{x}_3(t))$.}\label{LESOadrc-sim}
\end{figure}

\begin{figure}[ht]\centering
\subfigure[]
 {\includegraphics[width=4.1cm,height=4cm]{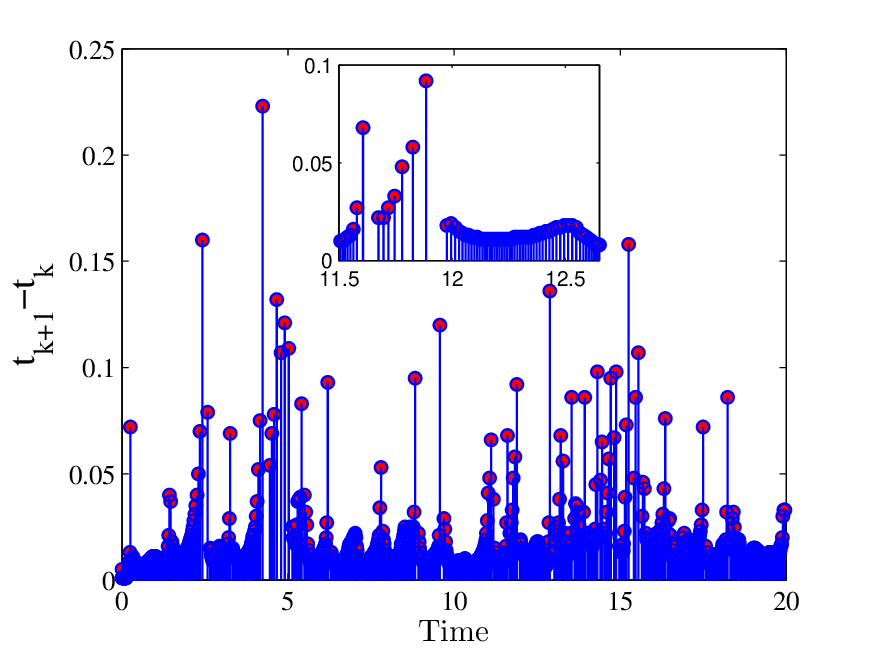}\label{Fig-ESOtime}}
\subfigure[]
 {\includegraphics[width=4.1cm,height=4cm]{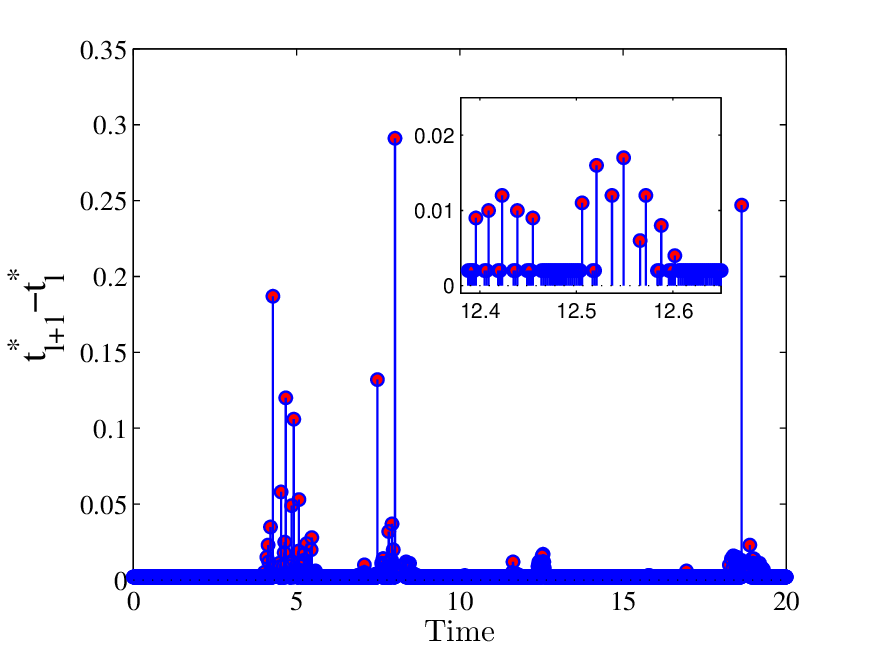}\label{Fig-Controlintertime}}
\caption{ Inter-execution times corresponding to ETMs \dref{ETMESOSimul} and \dref{ETMADRCSimul}.}\label{ESOADRC-inter}
\end{figure}

It can be observed from Figure \ref{LESOadrc-sim} that the estimate effect for both state $(x_{1}(t),x_{2}(t))$
and random total disturbance $x_{3}(t)$  and the stabilizing effect for $(x_{1}(t),x_{2}(t))$ are satisfactory,
where the estimation effect for $x_{3}(t)$ is not as good as the one for state $(x_{1}(t),x_{2}(t))$.
These are consistent with the theoretical result presented in Theorem \ref{theorem34.1}.
The inter-execution times corresponding to
ETM \dref{ETMESOSimul} for output transmission and ETM \dref{ETMADRCSimul} for control signal update
can be seen from Figure \ref{Fig-ESOtime} and Figure \ref{Fig-Controlintertime}, respectively,
whose respective number of execution times during [0, 20] is 1290 and 8554.

\section{Concluding remarks}\label{sectionV}
In this paper, event-triggered ADRC has been first addressed for a class of
uncertain random nonlinear systems.
The controlled systems with lower triangular  structure are subject to
nonlinear unmodeled dynamics, bounded noise, and colored noise in large scale,
whose total effects are treated as a random total disturbance. An event-triggered
ESO has been designed for real-time estimation of the random total disturbance,
and then an event-triggered controller composed of an output-feedback controller and a compensator has been
designed for the output-feedback stabilization and disturbance rejection for the controlled systems.
Rigorous theoretical proofs have been given to obtain both the mean square practical convergence and almost surely practical one
of the ADRC's closed-loop under two respective event-triggering mechanisms, validated by some numerical simulations.
Potential interesting problems to be solve are the design and convergence analysis of
 periodic event-triggered ADRC for uncertain random nonlinear systems and the comparison
of its communication efficiency with the one of the strategy proposed in this paper.

\bibliographystyle{plain}        

\end{document}